\documentclass[11pt]{article}

\usepackage{amsfonts,amsmath}
\usepackage{latexsym}
\usepackage{color}

\renewcommand{\theequation}{\arabic{section}.\arabic{equation}}

\newtheorem{prop}{Proposition}[section]
\newtheorem{cor}[prop]{Corollary}

\newtheorem{lem}[prop]{Lemma}
\newtheorem{theo}[prop]{Theorem}
\newtheorem{rema}[prop]{Remark}
\newtheorem{ex}[prop]{Example}

\newcommand{\R}{\mathbb{R}}
\newcommand{\E}{\mathbb{E}}
\newcommand{\PP}{\mathbb{P}}

\newcommand{\N}{\mathbb{N}}

\newcommand{\Z}{\mathbb{Z}}

\newcommand{\al}{\alpha}
\newcommand{\la}{\lambda}

\newcommand{\ga}{\gamma}

\newcommand{\ka}{\kappa}

\newcommand{\vpi}{\varpi}
\newcommand{\si}{\sigma}

\newcommand{\te}{\theta}

\newcommand{\be}{\beta}
\newcommand{\ep}{\varepsilon}

\newcommand{\de}{\delta}
\newcommand{\De}{\Delta}
\newcommand{\om}{\omega}
\newcommand{\Om}{\Omega}
\newcommand{\ze}{\zeta}

\newcommand{\f}{\mathcal{F}}
\newcommand{\g}{\mathcal{G}}
\newcommand{\h}{\mathcal{H}}

\newcommand{\laa}{\mathcal{L}}
\newcommand{\m}{\mathcal{M}}
\newcommand{\n}{\mathcal{N}}

\newcommand{\va}{\mathcal{V}}

\newcommand{\wc}{\widehat{c}}

\newcommand{\WE}{\widetilde{\mathbb{E}}}

\newcommand{\Wc}{\widetilde{c}}

\newcommand{\Wsi}{\widetilde{\sigma}}

\newcommand{\Wde}{\widetilde{\delta}}

\newcommand{\Wb}{\widetilde{b}}

\newcommand{\Bb}{\overline{b}{}}

\newcommand{\BL}{\overline{L}{}}

\newcommand{\BV}{\overline{V}{}}

\newcommand{\Bc}{\overline{c}{}}

\newcommand{\Bv}{\overline{v}{}}

\newcommand{\Bh}{\overline{h}{}}

\newcommand{\Bal}{\overline{\alpha}{}}
\newcommand{\Bbe}{\overline{\beta}{}}

\newcommand{\Bva}{\overline{\va}{}}

\def\nib{\noindent\bf}
\def\nit{\noindent\it}
\def\ni{\noindent}

\newcommand{\rdn}{\sqrt{\De_n}}

\newcommand{\rdnn}{\frac{1}{\rdn}}

\newcommand{\sdt}{\sum_{i=1}^{[t/\De_n]}}

\newcommand{\sdknt}{\sum_{i=1}^{[t/\De_n]-k_n+1}}

\newcommand{\dd}{\De^n_i}

\newcommand{\vsc}{\vskip 5mm}

\newcommand{\vst}{\vskip 3mm}
\newcommand{\vsq}{\vskip 4mm}

\newcommand{\toop}{\stackrel{\PP}{\longrightarrow}}

\newcommand{\tolfs}{\stackrel{\laa_f-\mbox{\tiny s}}{\longrightarrow}}
\newcommand{\tolls}{\stackrel{\laa-\mbox{\tiny s}}{\Longrightarrow}}

\newcommand{\tols}{\stackrel{\laa-\mbox{\tiny s}}{\longrightarrow}}

\newcommand{\toucp}{\stackrel{\mbox{\tiny u.c.p.}}{\Longrightarrow}}

\newcommand{\proba}{(\Omega ,\f,(\f_t)_{t\geq0},\PP)}

\newcommand{\probt}{(\widetilde{\Omega},\widetilde{\f},
(\widetilde{\f}_t)_{t\geq0},\widetilde{\PP})}

\newcommand{\bee}{\begin{equation}}
\newcommand{\eee}{\end{equation}}
\newcommand{\bea}{\begin{eqnarray}}
\newcommand{\eea}{\end{eqnarray}}
\newcommand{\bean}{\begin{eqnarray*}}
\newcommand{\eean}{\end{eqnarray*}}

\newcommand{\qed}{$\hfill\Box$}

\begin{document}

\author{Jean Jacod\thanks{Institut de Math\'ematiques de Jussieu, 4 Place
Jussieu, 75 005 Paris, France (CNRS -- UMR 7586, and Universit\'e Pierre et
Marie Curie), Email: jean.jacod@upmc.fr}
\and
Mathieu Rosenbaum\thanks{Laboratoire de Probabilit\'es et Mod\`eles
Al\'eatoires, 4 Place Jussieu, 75 005 Paris, France (CNRS -- UMR 7599,
and Universit\'e Pierre et Marie Curie), Email: mathieu.rosenbaum@upmc.fr}}

\title{Estimation of volatility functionals: the case of a $\sqrt{n}$
window}

\date{\today}

\maketitle

\begin{abstract}
We consider a multidimensional It\^o semimartingale regularly sampled on
$[0,t]$ at high frequency $1/\Delta_n$, with $\Delta_n$ going to zero. The
goal of this paper is to provide an estimator
for the integral over $[0,t]$ of a given function of the volatility matrix,
with the optimal rate $1/\rdn$ and minimal asymptotic variance. To achieve this
we use spot volatility estimators based on observations within time intervals
of length $k_n\De_n$. In \cite{JR} this was done with $k_n\to\infty$ and $k_n
\rdn\to0$, and a central limit theorem was given after suitable de-biasing.
Here we do the same with the choice $k_n\asymp1/\rdn$. This results in a
smaller bias, although more difficult to eliminate.
\end{abstract}



\section{Introduction}\label{sec-I}
\setcounter{equation}{0}
\renewcommand{\theequation}{\thesection.\arabic{equation}}

Consider an It\^o semimartingale $X_t$, whose squared volatility
$c_t$ (a $d\times d$ matrices-valued process if $X$ is $d$-dimensional)
is itself another It\^o semimartingale.
The process $X$ is observed at discrete times $i\De_n$ for
$i=0,1,\cdots$, the time lag $\De_n$ being small (high-frequency setting)
and eventually going to $0$. The aim is to estimate integrated functionals
of the volatility, that is $\int_0^tg(c_s)\,ds$ for arbitrary (smooth
enough) functions $g$, on the basis of the observations at stage
$n$ and within the time interval $[0,t]$.

In \cite{JR}, to which we refer for detailed motivations for this problem,
we have exhibited estimators which are consistent, and
asymptotically optimal, in the sense that they asymptotically achieve
the best rate $1/\rdn$, and also the minimal asymptotic variance in the cases
where optimality is well-defined (namely, when $X$ is continuous and has
a Markov type structure, in the sense of \cite{CDG}). These estimators have
this rate and minimal asymptotic variance as soon as the jumps of $X$
are summable, plus some mild technical conditions.

The aim of this report is to complement \cite{JR} with another estimator,
of the same type, but using spot volatility estimators based on a different
window size. In this introduction we explain the differences between the
estimator in that paper and the one presented here.

For the sake of simplicity we consider the case when $X$ is continuous and
one-dimensional (the discontinuous and multi-dimensional case is considered
later), that is of the form
$$X_t=X_0+\int_0^tb_s\,ds+\int_0^t\si_s\,dW_s$$
and $c_t=\si_t^2$ is the squared volatility. Natural estimators for
$\int_0^tg(c_s)\,ds$ are
\bee\label{I-2}
V^n(g)_t=\De_n\sdknt g(\wc^n_i),\quad
\text{where}~~\wc^n_i=\frac1{k_n\De_n}\sum_{j=0}^{k_n-1}
(X_{(i+j)\De_n}-X_{(i+j-1)\De_n})^2
\eee
for an arbitrary sequence of integers such that $k_n\to\infty$ and
$k_n\De_n\to0$: One knows that $V^n(g)_t\toop V(g)_t$ (when $g$
is continuous and of polynomial growth).

The variables $\wc^n_i$ are spot volatility estimators, and according to
\cite{JP} we know that $\wc^n_{[t/\De_n]}$ estimate $c_t$, with a rate
depending on the ``window size'' $k_n$. The optimal rate $1/\De_n^{1/4}$ is
achieved by taking $k_n\asymp1/\rdn$. When $k_n$ is smaller, the rate is
$\sqrt{k_n}$ and the estimation
error is a purely ``statistical error''; when $k_n$ is bigger, the rate is
$1/\sqrt{k_n\De_n}$ and the estimation error is due to the variability of
the volatility process $c_t$ itself (its volatility and its jumps). With the
optimal choice $k_n\asymp1/\rdn$ the estimation error is a mixture of the
statistical error and the error due to the variability of $c_t$.

In \cite{JR} we have used a ``small'' window, that is $k_n\ll 1/\rdn$.
Somewhat surprisingly, this allows for optimality in the estimation of
$\int_0^tg(c_s)\,ds$ (rate $1/\rdn$ and minimal asymptotic variance). However,
the price to pay is the need of a
de-biasing term to be subtracted from $V^n(g)$, without which the rate
is smaller and no Central Limit Theorem is available.

Here, we consider the window size $k_n\asymp1/\rdn$. This leads to a
convergence rate $1/\rdn$ for $V^n(g)$ itself, and the limit is again
conditionally Gaussian with the ``minimal'' asymptotic variance, but with
a bias that depends on the volatility of the volatility $c_t$, and on
its jumps. It is however possible to subtract from $V^n(g)$  a de-biasing
term again, so that the limit becomes (conditionally) centered.

Section \ref{sec-SET} is devoted to presenting assumptions and results,
and all proofs are gathered in Section \ref{sec-P}. The reader is referred
to \cite{JR} for motivation and various comments and a detailed
discussion of optimality. However, in order to make this report readable,
we basically give the full proofs, even though a number of partial results
have already been proved in the above-mentioned paper, and with the exception
of a few well designated lemmas.

\section{The results}\label{sec-SET}
\setcounter{equation}{0}
\renewcommand{\theequation}{\thesection.\arabic{equation}}

\subsection{Setting and Assumptions}
The underlying process $X$ is $d$-dimensional, and observed at the
times $i\De_n$ for $i=0,1,\cdots$, within a fixed interval of interest
$[0,t]$. For any process we write $\dd Y=Y_{i\De_n}-Y_{(i-1)\De_n}$
for the increment over the $i$th observation interval. We assume that the
sequence $\De_n$ goes to $0$. The precise assumptions on $X$ are as follows:

First, $X$ is an It\^o semimartingale on a filtered
space $\proba$. It can be written in its Grigelionis form,
as follows, using a $d$-dimensional Brownian motion $W$ and a Poisson
random measure $\mu$ on
$\R_+\times E$, with $E$ is an auxiliary Polish space and with
the (non-random) intensity measure $\nu(dt,dz)=dt\otimes\la(dz)$ for
some $\si$-finite measure $\la$ on $E$:
\bee\label{S-1}
\begin{array}{lll}
X_t&=&X_0+\int_0^tb_s\,ds+\int_0^t\si_s\,dW_{s}+\int_0^t\int_E
\de(s,z)\,1_{\{\|\de(s,z)\|\leq1\}}\,(\mu-\nu)(ds,dz)\\&&\hskip5cm
+\int_0^t\int_E\de(s,z)\,1_{\{\|\de(s,z)\|>1\}}\,\mu(ds,dz).
\end{array}
\eee
This is a vector-type notation: the process $b_t$ is $\R^d$-valued optional,
the process $\si_t$ is $\R^d\otimes\R^d$-valued optional, and
$\de=\de(\om,t,z)$ is a predictable $\R^d$-valued function on
$\Om\times\R_+\times E$.

The spot volatility process $c_t=\si_t\si_t^*$ ($^*$ denotes transpose)
takes its values in the set $\m^+_d$ of all nonnegative symmetric
$d\times d$ matrices. We suppose that $c_t$ is again an It\^o semimartingale,
which can be written as
\bee\label{S-2}
\begin{array}{lll}
c_t&=&c_0+\int_0^t\Wb_s\,ds+\int_0^t\Wsi_s\,dW_{s}+\int_0^t\int_E
\Wde(s,z)\,1_{\{\|\Wde(s,z)\|\leq1\}}\,(\mu-\nu)(ds,dz)\\
&&\hskip5cm+\int_0^t\int_E
\Wde(s,z)\,1_{\{\|\Wde(s,z)\|>1\}}\,\mu(ds,dz).
\end{array}
\eee
with the same $W$ and $\mu$ as in (\ref{S-1}). This is indeed {\em not a
restriction}: if $X$ and $c$ are two It\^o semimartingales, we have a
representation as above for the pair $(X,c)$ and, if the dimension of $W$
exceeds the dimension of $X$ one can always add fictitious component to
$X$, arbitrarily set to $0$, so that the dimensions of $X$
and $W$ agree.

In (\ref{S-2}), $\Wb$ and $\Wsi$ are optional and $\Wde$ is as $\de$; moreover
$\Wb$ and $\Wde$ are $\R^{d^2}$-valued. Finally, we need the spot volatility
of the volatility and ``spot covariation'' of the
continuous martingale parts of $X$ and $c$, which are
$$\Wc_t^{ij,kl}=\sum_{m=1}^d\Wsi_t^{ij,m}\Wsi_t^{kl,m},
\qquad
\Wc'^{i,jk}_t=\sum_{l=1}^d\si_t^{il}\Wsi_t^{jk,l}.$$

The precise assumption on the coefficients are as follows, with $r$
a real in $[0,1)$:
\vsq

\nib Assumption (A'-$r$): \rm There are a sequence $(J_n)$ of nonnegative
bounded $\la$-integrable functions on $E$ and a sequence
$(\tau_{n})$ of stopping times increasing to $\infty$, such that
\bee\label{S-3}
t\leq\tau_n(\om)~~\Longrightarrow~~
\|\de(\om,t,z)\|^r\wedge1+\|\Wde(\om,t,z)\|^2\wedge1 \leq J_n(z)
\eee
Moreover the processes $b'_t=b_t-\int\de(t,z)\,1_{\{\|\de(t,z)\|\leq1\}}\,
\la(dz)$ (which is well defined), $\Wc_t$ and $\Wc'_t$ are
c\`adl\`ag or c\`agl\`ad, and the maps $t\mapsto\Wde(\om,t,z)$ are
c\`agl\`ad (recall that $\Wde$ should be predictable), as well as
the processes $\Wb_t+\int\Wde(t,z)\big(\ka(\|\Wde(t,z)\|)-
1_{\{\|\Wde(t,z)\|\leq1\}})\,\la(dz)$ for one (hence for all) continuous
function $\ka$ on $\R_+$ with compact support and equal to $1$ on a
neighborhood of $0$.\qed
\vsq

The bigger $r$ is, the weakest Assumption (A-$r$) is, and when
(A-$0$) holds the process $X$ has finitely many jumps on each finite
interval. The part of (A-$r$) concerning the jumps of $X$ implies that
$\sum_{s\leq t}\|\De X_s\|^r<\infty$ a.s.\ for all $t<\infty$, and it
is in fact ``almost'' implied by this property. Since $r<1$, this
implies $\sum_{s\leq t}\|\De X_s\|<\infty$ a.s.

\begin{rema}\label{RS-1} \rm (A'-$r$) above is basically the same as
Assumption (A-$r$) in \cite{JR}, albeit (slightly) stronger (hence its name):
some degree of regularity in time seems to be
needed for $\Wb,\Wc,\Wc',\Wde$ in the present case.
\end{rema}

\subsection{A First Central Limit Theorem.}
For defining the estimators of the spot volatility, we first
choose a sequence $k_n$ of integers which satisfies, as $n\to\infty$:
\bee\label{R-1}
k_n~\sim~\frac{\te}{\rdn},\qquad \te\in(0,\infty),
\eee
and a sequence $u_n$ in $(0,\infty]$. The
$\m_d^+$-valued variables $\Wc^n_i$ are defined, componentwise, as
\bee\label{R-3}
\wc^{n,lm}_i=\frac1{k_n\De_n}\sum_{j=0}^{k_n-1}\De^n_{i+j}
X^l\,\De^n_{i+j} X^m\,1_{\{\|\De^n_{i+j}X\|\leq u_n\}},
\eee
and they implicitly depend on $\De_n,k_n,u_n$.

One knows that $\wc^n_{[t/\De_n]}\toop c_t$ for any $t$, and there is an
associated Central Limit Theorem under (A-$2$), with rate $1/\De_n^{1/4}$:
the choice (\ref{R-1}) is optimal, in the sense that it allows us to
have the fastest possible rate by a balance between the involved
``statistical error'' which is of order $1/\sqrt{k_n}$, and the
variation of $c_t$ over the interval $[t,t+k_n\De_n]$, which
is of order $\sqrt{k_n\De_n}$ because $c_t$ is an It\^o semimartingale
(and even when it jumps), see \cite{APPS,JP}.

By Theorem 9.4.1 of \cite{JP}, one also knows that under (A'-$r$) and
if $u_n\asymp \De_n^\vpi$ for some $\vpi\in\big[\frac{p-1}{2p-r},\frac12\big)$
we have
\bee\label{R-4}
V(g)^n_t:=\De_n\sdknt g(\wc^n_i)~\toucp~V(g)_t
:=\int_0^tg(c_s)\,ds
\eee
(convergence in probability, uniform over each compact interval; by
convention $\sum_{i=a}^bv_i=0$ if $b<a$), as soon as the function
$g$ on $\m^+_d$ is continuous with $|g(x)|\leq K(1+\|x\|^p)$ for some
constants $K,p$. Actually, for this to hold we need much weaker assumptions
on $X$, but we do not need this below. Note also that when $X$ is
continuous, the truncation in
(\ref{R-3}) is useless: one may use (\ref{R-3}) with $u_n\equiv\infty$,
which reduces to (\ref{I-2}) in the one-dimensional case.

Now, we want to determine at which rate the convergence (\ref{R-4})
takes place. This amounts to proving an associated Central Limit Theorem.
For an appropriate choice of the truncation
levels, such a CLT is available for $V(g)^n$, with the rate $1/\rdn$, but the
limit exhibits a bias term. Below, $g$ is a smooth function on
$\m^+_d$, and the two first partial derivatives are denoted as
$\partial_{jk}g$ and $\partial^2_{jk,lm}g$, since any $x\in\m^+_d$ has
$d^2$ components $x^{jk}$. The family of all partial derivatives of
order $j$ is simply denoted as $\partial^jg$.

\begin{theo}\label{TR-1} Assume (A'-$r$) for some $r<1$. Let $g$ be a $C^3$
function on $\m^+_d$ such that
\bee\label{R-8}
\|\partial^jg(x)\|\leq K(1+\|x\|^{p-j}),\qquad j=0,1,2,3
\eee
for some constants $K>0,\,p\geq3$. Either suppose that $X$ is continuous
and $u_n/\De_n^\ep\to\infty$ for some $\ep<1/2$ (for example,
$u_n\equiv\infty$, so there is no truncation at all), or suppose that
\bee\label{R-9}
u_n\asymp\De_n^\vpi,\qquad\frac{2p-1}{2(2p-r)}\leq\vpi<\frac12.
\eee
Then we have the finite-dimensional (in time) stable convergence in law
\bee\label{R-10}
\rdnn\,(V(g)^n_t-V(g)_t)~\tolfs~A^1_t+A^2_t+A^3_t+A^4_t+Z_t,
\eee
where $Z$ is a process defined on an extension
$\probt$ of $\proba$, which conditionally on $\f$ is a continuous
centered Gaussian martingale with variance
\bee\label{R-11}
\WE\big((Z_t)^2\mid\f\big)~=~\sum_{j,k,l,m=1}^d\int_0^t
\partial_{jk} g(c_s)\,\partial_{lm} g(c_s)\,
\big(c_s^{jl}c_s^{km}+c_s^{jm}c_s^{kl}\big)\,ds,
\eee
and where, with the notation
\bee\label{R-132}
G(x,y)=\int_0^1\big(g(x+wy)-wg(x+y)-(1-w)g(x)\big)\,dw,
\eee
we have
\bee\label{R-13}
\begin{array}{l}
A^1_t=-\frac{\te}2\,\big(g(c_0)+g(c_t)\big)\\
A^2_t=\frac1{2\te}\,\sum\limits_{j,k,l,m=1}^d
\int_0^t\partial^2_{jk,lm} g(c_s)
\,\big(c^{jl}_sc^{km}_s+c^{jm}_sc^{kl}_s\big)\,ds\\
A^3_t=-\frac{\te}{12}\,\sum\limits_{j,k,l,m=1}^d
\int_0^t\partial^2_{jk,lm} g(c_s)\,\Wc_s^{jk,lm}\,ds\\
A^4_t=\te\sum\limits_{s\leq t}G(c_{s-},\De c_s).
\end{array}
\eee
\end{theo}

Note that $|G(x,y)|\leq K(1+\|x\|)^p\,\|y\|^2$, so the sum defining
$A_t^4$ is absolutely convergent, and vanishes when $c_t$ is continuous.

\begin{rema}\label{RR-01} \rm  The bias has four parts:

1) The first one is due to a border effect: the formula giving $V^n(g)_t$
contains $[t/\De_n]-k_n+1$ summands only, whereas the natural (unfeasible)
approximation $\De_n\sdt g(c_{(i-1)\De_n})$ contains $[t/\De_n]$ summands. The
sum of the lacking $k_n$ summands is of order of magnitude $(k_n-1)\De_n$,
which goes to $0$ and thus does not impair
consistency, but it creates an obvious bias after normalization
by $1/\rdn$. Removing this source of bias is straightforward: since
$g(c_s)$ is ``under-represented'' when $s$ is close to $0$ or
to $t$, we add to $V^n(g)_t$ the variable
\bee\label{R-6}
\frac{(k_n-1)\De_n}2\,\big(g(\wc^n_1)+g(\wc^n_{[t/\De_n]-k_n+1}\big).
\eee
Of course, other weighted averages of $g(\wc^n_i)$ for $i$ close to $0$
or to $[t/\De_n]-k_n+1$ would be possible.

2)  The second part $A^2$ is continuous in time and is present even
for the toy model $X_t=\sqrt{c}\,W_t$ with $c$ a constant and
$\De_n=\frac1n$ and $T=1$. In this simple case it can be interpreted as
follows: instead of taking the ``optimal'' $g(\wc_n)$ for estimating $g(c)$,
with $\wc_n=\sum_{i=1}^n(\dd X)^2$, one takes $\frac1n\sum_{i=1}^n g(\wc^n_i)$
with $\wc^n_i$ a ``local'' estimator of $c$. This adds a statistical error
which results in a bias. Note that, even in the general case, this bias
would disappear, were we taking in (\ref{R-1}) the (forbidden) value
$\te=\infty$ (with still $k_n\De_n\to0$, at the expense of a slower rate
of convergence.

3) The third and fourth parts $A^3$ and $A^4$ are respectively
continuous and purely discontinuous, and due to the continuous part and
to the jumps of the volatility process $c_t$ itself. These two biases
disappear if we take $\te=0$ in (\ref{R-1}) (with still $k_n\to\infty$),
again a forbidden value,
and again at the expense of a slower rate of convergence.

\ni The only test function $g$ for which the last three biases disappear
is the identity $g(x)=x$. This is because, in this case, and up to the
border terms, $V(g)^n_t$ is nothing but the realized quadratic
variation itself and the spot estimators $\wc^n_i$ actually merge together
and disappear as such.
\end{rema}

\begin{rema}\label{RR-02} \rm Observe that
(\ref{R-9}) implies $r<1$. This restriction is not a surprise, since
one needs $r\leq1$ in order to estimate the integrated
volatility by the (truncated) realized volatility, with a
rate of convergence $1/\rdn$. When $r=1$ it is likely that the CLT
still holds fore an appropriate choice of the sequence $u_n$, and with
another additional bias, see e.g. \cite{V10} for a slightly different
context. Here we let this borderline case aside.
\end{rema}

\subsection{Estimation of the Bias.}

Now we proceed to ``remove'' the bias, which means subtracting
consistent estimators for the bias from $V'(g)^n_t$. As written before,
we have
\bee\label{R-601}
A^{n,1}_t=-\frac{k_n\rdn}2\,\big(g(\wc^n_1)
+g(\wc^n_{[t/\De_n]-k_n+1}\big)~\toop~A^1_t
\eee
(this comes from $\wc^n_1\toop c_0$ and $\wc^n_{[t/\De_n]-k_n+1}\toop
c_{t-}$, plus $c_{t-}=c_t$ a.s.).
Next, observe that $A^2=\frac1{\te}\,V(h)$ for the test
function $h$ defined on $\m^+_d$ by
\bee\label{R-100}
h(x)=\frac1{2}\,\sum_{j,k,l,m=1}^d\partial^2_{jk,lm} g(x)
\,\big(x^{jl}x^{km}+x^{jm}x^{kl}\big).
\eee
Therefore
\bee\label{R-101}
A^{n,2}_t=\frac1{k_n\rdn}\,V(h)^n_t~\toop~A^2_t.
\eee

The term $A_t^3$ involves the volatility of the volatility, for which
estimators have been provided in the one-dimensional case by M. Vetter in
\cite{V1}; namely, if $d=1$ and under suitable technical assumptions
(slightly stronger than here), {\em plus}\
the continuity of $X_t$ and $c_t$, he proves that
$$\frac3{2k_n}\sum_{i=1}^{[t/\De_n]-2k_n+1}
(\wc^n_{i+k_n}-\wc^n_i)^2$$
converges to $\int_0^t\big(\Wc_s+\frac6{\te^2}\,(c_s)^2\big)\,ds$. Of course, we need to
modify this estimator here, in order to include the function $\partial^2g$ in
the limit and account for the possibilities of having $d\geq2$ and
having jumps in $X$. We propose to take
\bee\label{R-15}
A^{n,3}_t=-\frac{\rdn}8\,\sum_{i=1}^{[t/\De_n]-2k_n+1}
\sum_{j,k,l,m=1}^d\partial^2_{jk,lm}g(\wc^n_i)\,
(\wc^{n,jk}_{i+k_n}-\wc^{n,jk}_i)\,(\wc^{n,lm}_{i+k_n}-\wc^{n,lm}_i).
\eee
When $X$ and $c$ are continuous one may expect the convergence to
$A_t^3-\frac12\,A_t^2$ (observe that
$\frac{\rdn}4\sim\frac3{2k_n}\,\frac{\te}{12}$), and one may expect the
same when $X$ jumps and $c$ is still continuous, because in (\ref{R-3}) the
truncation basically eliminates the jumps of $X$. In contrast,
when $c$ jumps, the limit should rather be related to the ``full''
quadratic variation of $c$, and indeed we have:

\begin{theo}\label{TR-2} Under the assumptions of Theorem \ref{TR-1},
for all $t\geq0$ we have
\bee\label{R-17}
A^{n,3}_t~\toop~-\frac12\,A^2_t+A^3_t+A'^4_t,
\eee
where
\bee\label{R-18}
A'^4_t=\te\sum_{s\leq t}G'(c_{s-},\De c_s)
\eee
and
\bee\label{R-180}
G'(x,y)=-\frac{1}8\sum_{j,k,l,m}\int_0^1
\big(\partial^2_{jk,lm}\,g(x)+
\partial^2_{jk,lm}\,g(x+(1-w)y)\big)\,w^2\,y^{jk}\,y^{lm}\,dw.
\eee
\end{theo}

At this stage, it remains to find consistent estimators for
$A^4_t-A'^4_t$, which has the form
$$A^4_t-A'^4_t=\te\sum_{s\leq t}G''(c_{s-},\De c_s),\quad
\text{where}~~ G''=G-G'.$$
More generally, we aim at
estimating
\bee\label{R-19}
\va(F)_t=\sum_{s\leq t}F(c_{s-},\De c_s),
\eee
at least when the function $F$ on $\m^+_d\times\m_d$, where $\m_d$
is the set of all $d\times d$ matrices, is $C^1$ and $|F(x,y)|\leq K\|y\|^2$
uniformly in $x$ within any compact set, as is the function $G''$ above.

The solution to this problem is not as simple as it might appear at first
glance. We first truncate from below, taking any sequence
$u'_n$ of truncation levels satisfying
\bee\label{R-20}
u'_n\to0,\qquad \frac{u'_n}{\De_n^{\vpi'}}\to\infty\quad\text{for
some}~~\vpi'\in\big(0,\frac18\big)
\eee
Second, we resort on the following trick. Since $\wc^n_i$ is
``close'' to the average of $c_t$ over the interval $(i\De_n,(i+k_n)\De_n]$,
we (somehow wrongly) pretend that, for all $j$:
\bee\label{R-21}
\begin{array}{c}
\exists s\in((j-1)k_n\De_n,jk_n\De_n]~~\text{with}~~\|\De c_s\|>u'_n~
\Leftrightarrow~\|\wc^n_{jk_n}-\wc^n_{(j-2)k_n}\|>u'_n\\
\De c_s\sim\wc^n_{jk_n}-\wc^n_{(j-2)k_n},\quad
\|\wc^n_{(j-1)k_n}-\wc^n_{(j-3)k_n}\|\bigvee\|\wc^n_{(j+1)k_n}-\wc^n_{(j-1)k_n}\|
<\|\wc^n_{jk_n}-\wc^n_{(j-2)k_n}\|.
\end{array}
\eee
The condition (\ref{R-20}) implies that for $n$ large enough there is
at most one jump of size bigger than $u'_n$ in each interval
$(i-1)\De_n,(i-1+k_n)\De_n]$ within $[0,t]$, and no two consecutive
intervals of this form contain such jumps. Despite this, the
statement above is of course not true, the main reason being that
$\wc^n_i$ and $c^n_i$ do not exactly agree. However it is ``true enough''
to allow for the next estimators to be consistent for $\va(F)_t$:
\bee\label{R-22}
\begin{array}{c}
\va(F)^n_t=\sum_{j=3}^{[t/k_n\De_n]-3}
F(\wc^n_{(j-3)k_n+1},\de^n_j\wc)\,1_{\{\|\de^n_{j-1}\wc\|\vee
\|\de^n_{j+1}\wc\|\vee u'_n<\|\de^n_j\wc\|\}},\\
\text{where}~~\de^n_j\wc=\wc^n_{jk_n+1}-\wc^n_{(j-2)k_n+1}.
\end{array}
\eee

Since this is a sum of approximately $[t/k_n\De_n]$ terms, the rate
of convergence of $\va(F)^n_t$ toward $\va(F)_t$ is law, probably
$1/\De_n^{1/4}$ only. However, here we are looking for consistent
estimators, and the rate is not of concern to us. Note that, again,
the upper limit in the sum above is chosen in such a way that
$\va(F)^n_t$ is computable on the basis of the observations within
the interval $[0,t]$.

\begin{theo}\label{TR-3} Assume all hypotheses of Theorem \ref{TR-1},
and let $F$ be a continuous function on $\R_+\times\R$ satisfying, with the
same $p\geq3$ as in (\ref{R-9}),
\bee\label{R-14}
|F(x,y)|\leq K(1+\|x\|+\|y\|)^{p-2}\,\|y\|^2.
\eee
Then for all $t\geq0$ we have
\bee\label{R-24}
\va(F)^n_t~\toop~\va(F)_t.
\eee
\end{theo}

\subsection{An Unbiased Central Limit Theorem.}

At this stage, we can set, with the notation (\ref{R-100}), (\ref{R-101}),
(\ref{R-15}) and (\ref{R-22}), and also (\ref{R-13}) and (\ref{R-180})
for $G$ and $G'$:
\bee\label{R-16}
\BV(g)^n_t=V(g)^n_t
+\frac{k_n\De_n}2\,\big(g(\wc^n_1)+g(\wc^n_{[t/\De_n]-k_n+1}\big)
-\rdn\,\Big(\frac32\,A^{n,2}_t
+A^{n,3}_t\Big)-k_n\De_n\,\va(G-G')^n_t.
\eee

We then have the following, which is a straightforward consequence
of the three previous theorems and of $k_n\rdn\to\te$, plus \eqref{R-601}
and \eqref{R-101} and the fact that
the function $G-G'$ satisfies (\ref{R-14}) when $g$ satisfies (\ref{R-8}):

\begin{theo}\label{TR-4} Under the assumptions of Theorem \ref{TR-1}, and
with $Z$ as in this theorem,
for all $t\geq0$ we have the finite-dimensional stable convergence in law
\bee\label{R-23}
\rdnn\,(\BV(g)^n_t-V(g)_t)~\tolfs~Z_t.
\eee
\end{theo}

Note that $\te$ no longer explicitly appears in this statement, so one can
replace (\ref{R-1}) by the weaker statement
\bee\label{R-29}
k_n~\asymp~\rdnn
\eee
(this is easily seen by taking subsequences $n_l$ such that
$k_{n_l}\,\sqrt{\De_{n_l}}$ converge to an arbitrary limit in $(0,\infty)$).

It is simple to make this CLT ``feasible'', that is usable in practice
for determining a confidence interval for $V(g)_t$ at any time $t>0$.
Indeed, we can define the following function on $\m^+_d$:
\bee\label{R-102}
\Bh(x)=\sum_{j,k,l,m=1}^d\partial_{jk}\,g(x)\,
\partial_{lm}\,g(x)\,\big(x^{jl}x^{km}+x^{jm}x^{kl}\big).
\eee
We then have $V(\Bh)^n\toucp V(\Bh)$, whereas $V(\Bh)_t$ is the right hand
side of (\ref{R-10}). Then we readily deduce:

\begin{cor}\label{CR-1} Under the assumptions of the previous theorem,
for any $t>0$ we have the following stable convergence in law,
where $Y$ is an $\n(0,1)$ variable:
\bee\label{R-25}
\frac{\BV(g)^n_t-V(g)_t}{\sqrt{\De_n\,V(\Bh)^n_t}}~\tols Y,\quad
\text{\rm in restriction to the set $\{V(\Bh)_t>0\}$},
\eee
\end{cor}

Finally, let us mention that the estimators $\BV(g)^n_t$ enjoy excatly
the same asymptotic efficiency properties as the estimators in
\cite{JR}, and we refer to this paper for a discussion of this topic.

\begin{ex}[Quarticity] \rm Suppose $d=1$ and take $g(x)=x^2$, so
we want tho estimate the quarticity $\int_0^tc_s^2\,ds$. In this case
we have
$$h(x)=2x^2,\qquad G(x,y)-G'(x,y)=-\frac{y^2}6.$$
Then the ``optimal'' estimator for the quarticity is
$$\De_n\big(1-\frac3{k_n}\big)\sdknt(\wc^n_i)^2
+\frac{\De_n}4\sum_{i=1}^{[t/\De_n]-2k_n+1}(\wc^n_{i+k_n}-\wc^n_i)^2
+\frac{(k_n-1)\De_n}2\,\big((\wc^n_1)^2+(\wc^n_{[t/\De_n]-k_n+1})^2\big).$$

The asymptotic variance is $8\int_0^tc_s^4\,ds$, to be compared with the
asymptotic variance of the more usual estimators
$\frac1{3\De_n}\sdt(\dd X)^4$, which is $\frac{32}3\int_0^tc_s^4\,ds$.
\end{ex}

\section{Proofs}\label{sec-P}
\setcounter{equation}{0}
\renewcommand{\theequation}{\thesection.\arabic{equation}}

\subsection{Preliminaries.}

According to the localization lemma 4.4.9 of \cite{JP} (for the assumption (K)
in that lemma), it is enough to show all four Theorems \ref{TR-1}, \ref{TR-2},
\ref{TR-3}, \ref{TR-4} under the following stronger assumption:
\vsq

\nib Assumption (SA'-$r$): \rm We have (A'-$r$). Moreover we have, for
a $\la$-integrable function $J$ on $E$ and a constant $A$:
\bee\label{P-0}
\|b\|,~\|b'\|,~\|\Wb\|,~\|c\|,~\|\Wc\|,~\|\Wc'\|,~J~\leq~A,\quad
\|\de(\om,t,z)\|^r\leq J(z),\quad \|\Wde(\om,t,z)\|^2\leq J(z),
\eee

In the sequel we thus suppose that $X$ satisfies (SA'-$r$), and
also that \eqref{R-1} holds: these assumptions are typically not recalled.
Below, all constants are denoted by $K$, and they vary from line to line.
They may implicitly depend on the process $X$ (usually through $A$ in
\eqref{P-0}). When they depend on an additional parameter $p$, we write $K_p$.

We will usually replace
the discontinuous process $X$ by the continuous process
\bee\label{P-3}
X'_t=\int_0^tb'_s\,ds+\int_0^t\si_s\,dW_{s},
\eee
connected with $X$ by $X_t=X_0+X'_t+\sum_{s\leq t}\De X_s$. Note that $b'$ is
bounded, and without loss of generality we will use below its c\`adl\`ag
version. Note also that, since the jumps of $c$ are bounded, one can rewrite
(\ref{S-2}) as
\bee\label{P-1}
c_t=c_0+\int_0^t\Wb_s\,ds+\int_0^t\Wsi_s\,dW_{s}+\int_0^t\int_E
\Wde(s,z)\,(\mu-\nu)(ds,dz).
\eee
This amounts to replacing $\Wb$ in (\ref{S-2}) by
$\Wb_{t+}+\int_E\de(t+,z)(\ka(\|\Wde(t+,z)\|)-1_{\{\|\Wde(t+,z)\|
\leq1\}})\,\la(dz)$, where $\ka$ is a continuous function with compact
support, equal to $1$ on the set $[0,A]$. Note that the new process
$\Wb$ is bounded c\`adl\`ag.

With any process $Z$ we associate the variables
\bee\label{P-502}
\eta(Z)_{t,s}=
\sqrt{\E\big(\sup\nolimits_{v\in(t,t+s]}\,\|Z_{t+v}-Z_t\|^2
\mid\f_t\big)},
\eee
and we recall Lemma 4.2 of \cite{JR}:

\begin{lem}\label{LP-1} For all $t>0$, all bounded c\`adl\`ag processes
$Z$, and all sequences $v_n\geq0$ of reals tending to $0$, we have
$\De_n\E\big(\sdt\eta(Z)_{(i-1)\De_n,v_n}\big)\to0$, and for all
$0\leq v\leq s$ we have $\E(\eta(Z)^n_{t+v,s}\mid\f_t)\leq\eta(Z)_{t,s}$.
\end{lem}

\subsection{An Auxiliary Result on It\^o Semimartingales}\label{ssec-ARIS}

In this subsection we give some simple estimates for a $d$-dimensional
semimartingale
$$Y_t=\int_0^tb^Y_s\,ds+\int_0^t\si^Y_s\,dW_{s}+\int_0^t\int_E
\de^Y(s,z)\,(\mu-\nu)(ds,dz)$$
on some space $\proba$, which may be different from the one on which $X$ is
defined, as well as $W$ and $\mu$, but we still suppose that the
intensity measure $\nu$ is the same. Note that $Y_0=0$ here. We assume
that for some constant $A$ and function $J^Y$ we have,
with $c^Y=\si^Y\si^{Y,*}$:
\bee\label{P-501}
\|b^Y\|\leq A,~\|c^Y\|\leq A^2,\quad
\|\de^Y(\om,t,z)\|^2\leq J^Y(z)\leq A^2,\quad\int_EJ^Y(z)\,\la(dz)\leq A^2.
\eee
The compensator of the quadratic variation of $Y$ is of the form
$\int_0^t\Bc^Y_s\,ds$, where
$\Bc^Y_t=c^Y_t+\int_E\de^Y(t,z)\,\de^Y(t,z)^*\,\la(dz)$. Moreover, if
the process $c^Y$ is itself an It\^o semimartingale, the quadratic
covariation of the continuous martingale parts of $Y$ and $c^Y$ is
also of the form $\int_0^t\wc'^Y_s\,ds$ for some process $\wc'^Y$,
necessarily bounded if both $Y$ and $c^Y$ satisfy (\ref{P-501}) (and,
if $Y=X$, we have $c^Y=c$ and $\wc'^Y=\wc'$).

\begin{lem}\label{LP-5} Below we assume (\ref{P-501}), and the constant
$K$ only depends on $A$.

a) We have for $t\in[0,1]$:
\bee\label{P-503}
\begin{array}{l}
\big\|\E(Y_t\mid\f_0)-tb^Y_0\big\|\leq t\,\eta(b^Y)_{0,t}\leq Kt\\
\big|\E(Y^j_t\,Y^m_t\mid\f_0)-t\Bc^{Y,jm}_0\big|\leq
Kt(t+\sqrt{t}\,\eta(b^Y)_{0,t}+\eta(\Bc^Y)_{0,t})\leq Kt,
\end{array}
\eee
and if further $\|\E(\Bc^Y_t-\Bc^Y_0\mid\f_0)\|\leq A^2t$ for all $t$,
we also have
\bee\label{P-505}
\big|\E(Y^j_t\,Y^m_t\mid\f_0)-t\Bc^{Y,jm}_0\big|\leq
2\,t^{3/2}(2A^2\sqrt{t}+A\eta(b^Y)_{0,t})\leq Kt^{3/2}.
\eee

b) When $Y$ is continuous, and if
$\E(\|\Bc^Y_t-\Bc^Y_0\|^2\mid\f_0)\leq A^4t$ for all $t$, we have
\bee\label{P-506}
\big|\E\big(Y^j_t\,Y^k_t\,Y^l_t\,Y^m_t\,\mid\f_0\big)
-t^2(c^{Y,jk}_0c^{Y,lm}_0+c^{Y,jl}_0c^{Y,km}_0
+c^{Y,jm}_0c^{Y,kl}_0)\big|\leq Kt^{5/2}.
\eee

c) When $c^Y$ is a (possibly discontinuous) semimartingale satisfying the
same conditions (\ref{P-501}) as $Y$, and if $Y$ itself is continuous, we
have
\bee\label{P-507}
\big|\E\big((Y^j_t\,Y^k_t-tc^{Y,jk}_0)(c_t^{Y,lm}-c_0^{Y,lm})\mid\f_0\big)
\leq Kt^{3/2}(\sqrt{t}+\eta(\wc'^Y)_{0,t}).
\eee

\end{lem}

\nib Proof. \rm The first part of (\ref{P-503}) follows by taking
the $\f_0$-conditional expectation in the decomposition
$Y_t=M_t+tb^Y_0+\int_0^t(b^Y_s-b^Y_0)\,ds$, where $M$ is a $d$-dimensional
martingale with $M_0=0$. For the second part, we deduce from It\^o's
formula that $Y^jY^m$ is the sum of a martingale vanishing at $0$ and of
$$b^j_0\!\int_0^t\!Y_s^m\,ds+b^m_0\!\int_0^t\!Y_s^j\,ds
+\int_0^t\!Y_s^m(b_s^j-b_0^j)\,ds+\int_0^t\!Y_s^j(b_s^m-b_0^m)\,ds
+\Bc^{Y,jm}_0t+\int_0^t\!(\Bc^{Y,jm}_s-\Bc^{Y,jm}_0)\,ds.$$
Since $\E(\|Y_t\|\mid\f_0)\leq KA\sqrt{t}$, as in (\ref{P-251}), we deduce
the second part of (\ref{P-503}) and also (\ref{P-505})
by taking again the conditional
expectation and by using the Cauchy-Schwarz inequality and the first part.

\eqref{P-506} is a part of Lemma 5.1 of \cite{JR}. For \eqref{P-507}, we first
observe that $Y^j_tY^k_t-tc^{Y,jk}_0=B_t+M_t$ and
$c_t^{Y,lm}-c_0^{Y,lm}=B'_t+M'_t$, with $M$ and $M'$ martingales
($M$ is continuous). The processes $B$, $B'$,
$\langle M,M\rangle$, $\langle M',M'\rangle$ and $\langle M,M'\rangle$
are absolutely continuous, with densities $\Bb_s$, $\Bb'_s$,
$h_s$, $h'_s$ and $h''_s$ satisfying, by (\ref{P-501}) for $Y$
and $c^Y$:
$$|\Bb_s|\leq 2\|Y_s\|\,\|b^Y_s\|+\|c^Y_s-c^Y_0\|,\quad
|\Bb'_s|\leq K,\quad
|h_s|\leq K\|Y_s\|^2,\quad|h'_s|\leq K,$$
whereas $h''_s=Y^j_s\wc'^{Y,k,lm}+Y^k_s\wc'^{Y,j;lm}$.
As seen before, $\E(\|Y_t\|^q\mid\f_0)\leq K_qt^{q/2}$
for all $q$, and $\E(\|c^Y_t-c_0^Y\|^2\mid\f_0)\leq Kt$. This yields
$\E(B_t^2\mid\f_0)\leq Kt^3$ and $\E(M_t^2\mid\f_0)\leq Kt^2$. Since
$|B'_t|\leq Kt$ and $\E(M'^2_t\mid\f_0)\leq Kt$, we deduce that
the $\f_0$- conditional expectations of $B_tB'_t$ and $B_tM'_t$ and
$M_tB'_t$ are smaller than $Kt^2$.

Finally $\E(M_tM'_t\mid\f_0)=\E(\langle M,M'\rangle_t\mid\f_0)$, and
$\langle M,M'\rangle_t$ is the sum of
$\wc'^{Y,k,lm}_0\int_0^tY^j_s\,ds+
\int_0^tY^j_s(\wc'^{Y,k,lm}_s-\wc'^{Y,k,lm}_0)\,ds$ and a similar term
with $k$ and $j$ exchanged. Then using again
$\E(\|Y_t\|^2\mid\f_0)\leq Kt$, plus $\|\E(Y_t\mid\f_0)\|\leq Kt$ and
Cauchy-Schwarz inequality, we obtain that the above conditional
expectation is smaller than $K(t^2+t^{3/2}\eta(\wc'^Y)_t)$. This
completes the proof of (\ref{P-507}).\qed

\subsection{Some Estimates.}

\nib 1) \rm We begin with well known estimates for $X'$
and $c$, under (\ref{P-0}) and for $s,t\geq0$ and $q\geq0$:
\bee\label{P-251}
\begin{array}{ll}
\E\big(\sup_{w\in[0,s]}\,\|X'_{t+w}-X'_t\|^q\mid\f_t\big)\leq K_q\,s^{q/2},
\quad&\|\,\E(X'_{t+s}-X'_t\mid\f_s)\|\leq Ks\\
\E\big(\sup_{w\in[0,s]}\,\|c_{t+w}-c_t\|^q\mid\f_t\big)
\leq K_q\,s^{1\wedge(q/2)},&
\|\,\E(c_{t+s}-c_t\mid\f_s)\|\leq Ks.
\end{array}
\eee

Next, it is much easier (although unfeasible in practice) to
replace $\wc^n_i$ in (\ref{R-4}) by the estimators based on the process $X'$
given by (\ref{P-3}). Namely, we will replace $\wc^n_i$ by the following:
$$\wc'^n_i=\frac1{k_n\De_n}\sum_{j=0}^{k_n-1}
\De^n_{i+j}X'\,\De^n_{i+j}X'^*.$$
The difference between $\wc^n_i$ and $\wc'^n_i$ is estimated by the following
inequality, valid when $u_n\asymp\De_n^\vpi$ and $q\geq1$, and where $a_n$
denotes a sequence of numbers (depending on $u_n$), going to $0$ as
$n\to\infty$ (this is Equation 4.8 of \cite{JR}):
\bee\label{P-6}
\E\big(\|\wc^n_i-\wc'^n_i\|^q\big) \leq K_qa_n\,\De_n^{(2q-r)\vpi+1-q}.
\eee
\vsq

\nib 2) \rm The jumps of $c$ also potentially cause troubles. So we will
eliminate the ``big'' jumps as follows. For any $\rho>0$ we
consider the subset $E_\rho=\{z:J(z)>\rho\}$, which satisfies
$\la(E_\rho)<\infty$, and we denote by $\g^\rho$ the
$\si$-field generated by the variables $\mu([0,t]\times A)$, where $t\geq0$
and $A$ runs through all Borel subsets of $E_\rho$. The process
\bee\label{P-601}
N^\rho_t=\mu((0,t]\times E_\rho)
\eee
is a Poisson process and we let $S^\rho_1,S^\rho_2,\cdots$ be its successive
jump times, and $\Om_{n,t,\rho}$ be the set on which
$S^\rho_j\notin\{i\De_n:\,i\geq1\}$ for all $j\geq1$ such that $S^\rho_j<t$,
and $S^\rho_{j+1}>t\wedge S^\rho_j+(6k_n+1)\De_n$  for all $j\geq0$
(with the convention $S_0^\rho=0$; taking $6k_n$ here instead of the more
natural $k_n$ will be needed in the proof of Theorem \ref{TR-3}, and makes
no difference here). All these objects are
$\g^\rho$-measurable, and $\PP(\Om_{n,t,\rho})\to1$ as
$n\to\infty$, for all $t,\rho>0$.

We define the processes
$$\Wb(\rho)_t=\Wb_t-\int_{E_\rho}\Wde(t+,z)\,\la(dz),\qquad
\Bc(\rho)_t=\Wsi_t\,\Wsi_t^*
+\int_{(E_\rho)^c}\Wde(t+,z)\,\Wde(t+,z)^*\,\la(dz)$$
\bee\label{P-101}
\begin{array}{lll}
c(\rho)_t&=&c_t-\int_0^t\int_{E_\rho}\Wde(s,z)\,\mu(ds,dz)~=~c^{(1)}(\rho)_t
+c^{(2)}(\rho)_t,\quad\text{where}\\
&&c^{(1)}(\rho)_t=c_0+\int_0^t\Wb(\rho)_s\,ds+\int_0^t\Wsi_s\,dW_s\\
&&c^{(2)}(\rho)_t=
\int_0^t\int_{(E_\rho)^c}\Wde(t-,z)\,(\mu-\nu)(ds,dz),
\end{array}
\eee
so $\Bc(\rho)$, which is $\R^{d^2}\otimes\R^{d^2}$-valued, is the
c\`adl\`ag version of the density of the
predictable quadratic variation of $c(\rho)$. Moreover $\g^\rho=
\{\emptyset,\Om\}$ and $(\Wb(\rho),c(\rho))=(\Wb,c)$ when $\rho$
exceeds the bound of the function $J$. Note also that $\Wb(\rho)$ and
$\Bc(\rho)$ are c\`adl\`ag.

By Lemma 2.1.5 and Proposition 2.1.10 in \cite{JP} applied to each
components of $X'$ and $c^{(2)}(\rho)$,
plus the property $\|\Wb(\rho)\|\leq K/\rho$, for all
$t\geq0$, $s\in[0,1]$, $\rho\in(0,1]$, $q\geq2$, we have
\bee\label{P-11}
\begin{array}{l}
\E\big(\sup_{w\in[0,s]}\,\|X'_{t+w}-X'_t\|^q\mid\f_t\vee\g^\rho\big)
\leq K_q\,s^{q/2}\\
\|\,\E(X'_{t+s}-X'_t\mid\f_s\vee\g^\rho)\|+
\|\,\E(c(\rho)_{t+s}-c(\rho)_t\mid\f_s\vee\g^\rho)\|\leq Ks\\

\E\big(\sup_{w\in[0,s]}\,
\|c^{(2)}(\rho)_{t+w}-c^{(2)}(\rho)_t\|^q\mid\f_t\vee\g^\rho\big)
\leq K_q\,\phi_\rho\,(s+s^{q/2})\\
\E\big(\sup_{w\in[0,s]}\,
\|c(\rho)_{t+w}-c(\rho)_t\|^q\mid\f_t\vee\g^\rho\big)
\leq K_q\big(\phi_\rho\,s+s^{q/2}+\frac{s^q}{\rho^q}\big)
\leq K_{q,\rho}\,s.
\end{array}
\eee
where
$\phi_\rho=\int_{(E_\rho)^c}J(z)\,\la(dz)\to0$ as $\rho\to0$. Note also
that $\|\Wb(\rho)_t\|\leq K/\rho$.
\vsc

\nib 3) \rm For convenience, we put
\bee\label{P-2}
\begin{array}{lll}
b^n_i=b_{(i-1)\De_n},\quad& c^n_i=c_{(i-1)\De_n}\\
\Wb(\rho)^n_i=\Wb(\rho)_{(i-1)\De_n},\quad&
\Bc(\rho)^n_i=\Bc(\rho)_{(i-1)\De_n},\quad&
c(\rho)^n_i=c(\rho)_{(i-1)\De_n}\\
\f^n_i=\f_{(i-1)\De_n},\quad&\f^{n,\rho}_i=\f_i^n\vee\g^\rho.
\end{array}
\eee
All the above variables are $\f_i^{n,\rho}$-measurable. Recalling
\eqref{P-502}, and writing $\eta(Z,(\h_t))_{t,s}$ if we use the
filtration $(\h_t)$ instead of $(\f_t)$, we also set
$$\eta(\rho)^n_{i,j}=\max(\eta(Y,(\g^\rho\bigvee\f_t))_{(i-1)\De_n,j\De_n}:
~Y=b',\Wb(\rho),c,\Bc(\rho),\wc'\big),
\quad \eta(\rho)^n_i=\eta(\rho)^n_{i,i+2k_n}.$$
Therefore, Lemma \ref{LP-1} yields for all $t,\rho>0$ and $j,k$ such that
$j+k\leq2k_n$:
\bee\label{P-511}
\De_n\E\big(\sdt\eta(\rho)^n_i\big)\to0,\qquad
\E(\eta(\rho)^n_{i+j,k}\mid\f^{n,\rho}_i)\leq\eta(\rho)^n_i.
\eee

We still need some additional notation.
First, define $\g^\rho$-measurable (random) set of integers:
\bee\label{P-134}
L(n,\rho)=\{i=1,2,\cdots:\,N^\rho_{(i+2k_n)\De_n}
-N^\rho_{(i-1)\De_n}=0\}
\eee
(taking above, $2k_n$ instead of $k_n$, is necessary for the proof of
Theorem \ref{TR-2}). Observe that
\bee\label{P-135}
i\in L(n,\rho),~0\leq j\leq 2k_n+1
~\Rightarrow~c^n_{i+j}-c^n_i=c(\rho)^n_{i+j}-c(\rho)^n_i.
\eee

Second, we define the following $\R^d\otimes\R^d$-valued variables
\bee\label{P-20}
\begin{array}{l}
\al^n_i=\dd X'\,\dd X'^*-c^n_i\,\De_n\\
\be^n_i=\wc'^n_i-c^n_i=\frac1{k_n\De_n}\,\sum_{j=0}^{k_n-1}
\big(\al^n_{i+j}+(c^n_{i+j}-c^n_i)\De_n\big)\\
\ga^n_i=\wc'^n_{i+k_n}-\wc'^n_i=\be^n_{i+k_n}-\be^n_i+c^n_{i+k_n}-c^n_i.
\end{array}
\eee
\vst

\nib 4) \rm Now we proceed with estimates.
(\ref{P-11}) yields, for all $q\geq0$:
\bee\label{P-21}
\begin{array}{ll}
\E(\|\al^n_i\|^q\mid\f^{n,\rho}_i)\leq K_q\De_n^q,\qquad&
\|\E(\al^n_i\mid\f^{n,\rho}_i)\|\leq K\De_n^{3/2}\\
\E\big(\|\sum_{j=0}^{k_n-1}\al^n_{i+j}\|^q\mid\f^{n,\rho}_i\big)
\leq K_q\De_n^{3q/4},~~&\E(\|\wc'^n_i\|^q\mid\f^{n,\rho}_i)\leq K_q
\end{array}
\eee
the third inequality following from the two first one, plus Burkholder-Gundy
and H\"older inequalities, and the last inequality form the third one and
the boundedness of $c_t$. Moreover, since the set $\{i\in L(n,\rho)\}$ is
$\g^\rho$-measurable,
the last part of (\ref{P-11}), (\ref{P-135}), and H\"older's
inequality, readily yield
\bee\label{P-1401}
q\geq2,~i\in L(n,\rho)~~\Rightarrow~~
\E\big(\|\be^n_i\|^q\mid\f^{n,\rho}_i)\big|\leq K_q
\Big(\rdn\,\phi_\rho+\De_n^{q/4}+\frac{\De_n^{q/2}}{\rho^q}\Big).
\eee
\vst

\nib 5) \rm The previous estimates are not enough for us.
We will apply the estimates of Lemma
\ref{LP-5} with $Y_t=X'_{(i-1)\De_n+t}-X'_{(i-1)\De_n}$
for any given pair $n,i$, and with the filtration
$(\f_{(i-1)\De_n+t}\vee\g^\rho)_{t\geq0}$. We observe that on the set
$A(\rho,n,i)=\{\exists j\leq2k_n:i-j\in L(n,\rho)\}$, which is
$\g^\rho$-measurable,
and because of (\ref{P-135}), the process $c^Y$ coincide with
$c(\rho)_{(i-1)\De_n+t}-c(\rho)_{(i-1)\De_n}$ if $t\in[0,\De_n]$.
Then in restriction to this set, by (\ref{P-505}) and (\ref{P-506})
and by the definition of $\eta(\rho)^n_{i,1}$, we have
$$\!\!\begin{array}{l}
\big|\E(\dd X'^j\,\dd X'^m\mid\f^{n,\rho}_i)-c^{n,jm}_i\De_n\big|
\leq K_\rho\De_n^{3/2}(\rdn+\eta(\rho)^n_{i,1})\\
\big|\E\big(\dd X'^j\,\dd X'^k\,\dd X'^l\,\dd X'^m\,\mid\f^{n,\rho}_i\big)
-(c^{n,jk}_ic^{n,lm}_i+c^{n,jl}_ic^{n,km}_i
+c^{n,jm}_ic^{n,kl}_i)\De_n^2\big|\leq K_\rho\De_n^{5/2}
\end{array}$$
(the constant above depends on $\rho$, through the bound $K/\rho$ for the
drift of $c(\rho)$). Then a simple calculation gives us
\bee\label{P-22}
\left.\begin{array}{l}
\big\|\,\E(\al^n_i\mid\f^{n,\rho}_i)\big\|\leq K_\rho\De_n^{3/2}
(\rdn+\eta(\rho)^n_{i,1})\\
\big|\E\big(\al^{n,jk}_i\al_i^{n,lm}\mid\f^{n,\rho}_i)
-(c^{n,jl}_ic^{n,km}_i+c^{n,jm}_ic^{n,kl}_i)\De_n^2\big|\leq K_\rho\De_n^{5/2}
\end{array}\right\}\quad\text{on}~~A(\rho,n,i)
\eee

Next, we apply Lemma \ref{LP-5} to the
process $Y_t=c(\rho)_{(i-1)\De_n+t}-c(\rho)_{(i-1)\De_n}$ for any given
pair $n,i$, and with the filtration $(\f_{(i-1)\De_n+t}\vee\g^\rho)_{t\geq0}$.
We then deduce from (\ref{P-503}), plus again (\ref{P-135}), that
\bee\label{P-13}
\begin{array}{l}
\hskip4cm i\in L(n,\rho),~0\leq t\leq k_n\De_n\Rightarrow\\
\big|\E((c^{jk}_{(i-1)\De_n+t}-c^{jk}_{(i-1)\De_n})
(c^{lm}_{(i-1)\De_n+t}-c^{lm}_{(i-1)\De_n})
\mid\f^{n,\rho}_i)-t\Bc(\rho)^{n,jklm}_i\big|
\leq K_\rho t\,\eta(\rho)^n_{i,k_n}\\
\big\|\E(c_{(i-1)\De_n+t}-c_{(i-1)\De_n}\mid\f^{n,\rho}_i)
-t\Wb(\rho)^n_i\,\big\|
\leq K_\rho t\,\eta(\rho)^n_{i,k_n}\leq K_p\,t.
\end{array}
\eee
Moreover, the Cauchy-Schwarz inequality and (\ref{P-21}) on the one
hand, and (\ref{P-507}) applied with the process
$Y_t=X'_{(i-1)\De_n+t}-X'_{(i-1)\De_n}$ on the other hand, give us
\bee\label{P-1301}
i\in L(n,\rho)~\Rightarrow\left\{
\begin{array}{l}
\big|\E\big(\al^{n,kl}_i\,\dd \Wb(\rho)^{ms}\mid\f^{n,\rho}_i\big)\big|
\leq K\De_n\eta(\rho)^n_{i,1}\\
\big|\E\big(\al^{n,kl}_i\,\dd c^{ms}\mid\f^{n,\rho}_i\big)
\big|\leq K_\rho\De_n^{3/2}(\rdn+\eta(\rho)^n_{i,1}).
\end{array}\right.
\eee
\vst

\nib 6) \rm We now proceed to estimates on $\be^n_i$:

\begin{lem}\label{LP-8} We have on the set where $i$ belongs to $L(n,\rho)$:
$$\begin{array}{l}
\big|\E(\be^{n,jk}_i\,\be^{n,lm}_i
\mid\f^{n,\rho}_i)-\frac1{k_n}\,(c^{n,jl}_ic^{n,km}_i+c^{n,jm}_ic^{n,kl}_i)
-\frac{k_n\De_n}3\,\Bc(\rho)_i^{n,jklm}\big|\\
\hskip9cm \leq K_\rho\rdn\,(\De_n^{1/4}+\eta(\rho)^n_i)\\
\big|\E(\be^{n,jk}_i(c^{n,lm}_{i+k_n}-c^{n,lm}_i)
\mid\f^{n,\rho}_i)-\frac{k_n\De_n}2\,\Bc(\rho)_i^{n,jklm}\big|
\leq K_\rho\rdn\,(\rdn+\eta(\rho)^n_i).
\end{array}$$
\end{lem}

\nib Proof. \rm We set $\ze^n_{i,j}=\al^n_{i+j}+(c^n_{i+j}-c_i^n)\De_n$ and
write $\be^{n,jk}_i\be_i^{n,lm}$ as
\bee\label{P-2201}
\frac1{k_n^2\De_n^2}
\sum_{u=0}^{k_n-1}\ze_{i,u}^{n,jk}\ze^{n,lm}_{i,u}
+\frac1{k_n^2\De_n^2}\sum_{u=0}^{k_n-2}\sum_{v=u+1}^{k_n-1}
\ze^{n,jk}_{i,u}\ze^{n,lm}_{i,v}
+\frac1{k_n^2\De_n^2}\sum_{u=0}^{k_n-2}\sum_{v=u+1}^{k_n-1}
\ze^{n,lm}_{i,u}\ze^{n,jk}_{i,v}.
\eee
For the estimates below, we implicitly assume $i\in L(n,\rho)$ and
$u,v\in\{0,\cdots,k_n-1\}$.

First, we deduce from \eqref{P-22} and \eqref{P-13}, plus \eqref{P-1301} and
successive conditioning, that
\bee\label{P-2208}
\big|\E(\ze_{i,u}^{n,jk}\ze^{n,lm}_{i,u}\mid\f^{n,\rho}_i)-
(c^{n,jl}_ic^{n,km}_i+c^{n,jm}_ic^{n,kl}_i)\De_n^2
\big|\leq K\De_n^{5/2}.
\eee
Second, if $u<v$, the same type of arguments and the boundedness of
$\Wb(\rho))_t$ and $c_t$ yield
$$\begin{array}{l}
|\E(\ze^{n,jk}_{i,v}\mid\f^{n,\rho}_{i+u+1})
-(c^{n,jk}_{i+u+1}-c^{n,jk}_i)\De_n-\Wb(\rho)^{n,jk}_{i+u+1}\De_n^2(v-u-1)|\\
\hskip7cm\leq K\De_n^{3/2}(k_n\rdn+\eta(\rho)^n_{i+v,1})\\
|\E(\al^{n,lm}_{i+u}\,(c^{n,jk}_{i+u+1}-c^{n,jk}_{i+u})\mid\f^{n,\rho}_{i+u})|
\leq K_\rho\De_n^{3/2}(\rdn+\eta(\rho)^n_{i+u,1})\\
|\E(\al^{n,lm}_{i+u}\,(c^{n,jk}_{i+u}-c^{n,jk}_i)\mid\f^{n,\rho}_{i+u})|
\leq K\De_n^{3/2}(\rdn+\eta^n_{i+u,1})\\
|\E(\al^{n,lm}_{i+u}\,(\Wb(\rho)^{n,jk}_{i+u+1}-\Wb(\rho)^{n,jk}_{i+u})
\mid\f^{n,\rho}_{i+u})|
\leq K_\rho\De_n^{3/2}(\rdn+\eta(\rho)^n_{i+u,1})\\
|\E(\al^{n,lm}_{i+u}\,\Wb(\rho)^{n,jk}_{i+u}\mid\f^{n,\rho}_{i+u})|
\leq K_\rho\De_n^{3/2}(\rdn+\eta(\rho)^n_{i+u,1})\\
|\E((c^{n,lm}_{i+u}-c^{n,lm}_i)\,(c^{n,jk}_{i+u+1}-c^{n,jk}_i)
\mid\f^{n,\rho}_{i})-\Bc(\rho)_i^{n,jklm}\De_nu|
\leq K_\rho\De_n\,\eta(\rho)^n_i\\
|\E((c^{n,lm}_{i+u}-c^{n,lm}_i)\,\Wb(\rho)^{n,jk}_{i+u+1}
\mid\f^{n,\rho}_{i})|\leq K_\rho\De_n^{1/4}.
\end{array}$$
Since $\sum_{u=0}^{k_n-2}\sum_{v=u+1}^{k_n-1}u=k_n^3/6+$ O$(k_n^2)$,
we easily deduce that the $\f^{n,\rho}_{i}$-conditional expectation of
the last term in \eqref{P-2201} is $\frac16\,\Bc(\rho)_i^{n,jklm}k_n\De_n$,
up to a remainder term which is O$(\rdn\,(\De_n^{1/4}+\eta(\rho)^n_i))$,
and the same is obviously true of the second term. The first claim of the
lemma readily follows from this and \eqref{P-2201} and \eqref{P-2208}.

The proof of the second claim is similar. Indeed, we have
$$\be_i^{n,jk}(c_{i+k_n}^{n,lm}-c_i^{n,lm})=
\frac1{k_n\De_n}
\sum_{u=0}^{k_n-1}\big(\al^{n,jk}_{i,u}+(c^{n,jk}_{i+u}-
c^{n,jk}_i)\De_n\big)\,\big(c^{n,lm}_{i+k_n}-c^{n,lm}_i)\big)$$
and
$$\big|\E(c^{n,lm}_{i+k_n}-c^{n,lm}_i\mid\f^{n,\rho}_{i+u+1})-
c^{n,lm}_{i+u+1}-c^{n,lm}_i-\Wb(\rho)^{n,lm}_{i+u+1}\De_n(k_n-u-1)
\big|\leq K\De_n\eta(\rho)^n_{i+u+1,k_n-u}.$$
Using the previous estimates, we conclude as for the first claim.\qed
\vsc

Finally, we deduce the following two estimates on the variables $\ga^n_i$ of
(\ref{P-20}), for any $q\geq2$:
\bee\label{P-26}
i\in L(n,\rho)~\Rightarrow\left\{\begin{array}{l}
\big|\E\big(\ga^{n,jk}_i\,\ga^{n,lm}_i\mid\f^{n,\rho}_i)
-\frac2{k_n}\,(c^{n,jl}_ic^{n,km}_i+c^{n,jm}_ic^{n,kl}_i)\\
\hskip2cm
-\frac{2k_n\De_n}3\,\Bc(\rho)^{n,jklm}_i\big|
\leq K_\rho\,\rdn\,\big(\De_n^{1/8}+\eta(\rho)^n_i\big)\\
\E(\|\ga^n_i\|^q\mid\f^{n,\rho}_i)\leq K_q
\big(\rdn\,\phi_\rho+\De_n^{q/4}+\frac{\De_n^{q/2}}{\rho^q}\big).
\end{array}\right.
\eee
To see that the first claim holds, one expands the product
$\ga^{n,jk}_i\,\ga^{n,lm}_i$ and use successive conditioning, the
Cauchy-Schwarz inequality and \eqref{P-11}, \eqref{P-135} and(\ref{P-13}),
and Lemma \ref{LP-8}; the contributing terms are
$$\begin{array}{c}
\be^{n,jk}_i\,\be^{n,lm}_i+\be^{n,jk}_{i+k_n}\,\be^{n,lm}_{i+k_n}
+(c^{n,jk}_{i+k_n}-c^{n,jk}_i)(c^{n,lm}_{i+k_n}-c^{n,lm}_i)\\
-\be^{n,jk}_i(c^{n,lm}_{i+k_n}-c^{n,lm}_i)
-\be^{n,lm}_i(c^{n,jk}_{i+k_n}-c^{n,jk}_i).
\end{array}$$
For the second claim we use (\ref{P-11}), \eqref{P-135} and
(\ref{P-1401}), and it holds for all $q\geq2$.

\subsection{The Behavior of Some Functionals of $c(\rho)$.}

For $\rho>0$ we set
\bee\label{P-30}
\begin{array}{c}
U(\rho)^{n}_t=\sum_{j=3}^{[t/k_n\De_n]-3}\|\mu(\rho)^{n}_j\|^2\,
1_{\{\|\mu(\rho)^n_j\|>u'_n/4\}},\quad\text{where}\\
\mu(\rho)^n_j=\frac1{k_n}\,\sum_{w=0}^{k_n-1}(c(\rho)^n_{jk_n+w}
-c(\rho)^n_{(j-2)k_n+w}).
\end{array}
\eee
The aim of this subsection is to prove the following lemma:

\begin{lem}\label{LP-2} Under (SA'-$r$) and (\ref{R-20}) we have
$$\lim_{\rho\to0}~\limsup_{n\to\infty}~\E\big(U(\rho)^n_t\big)~=~0.$$
\end{lem}

Assumption (SA'-$r$) is of course not fully used. What is
needed is the assumptions concerning the process $c_t$ only.
\vsq

\nib Proof. \rm With the notation (\ref{P-101}), and for $l=1,2$ we define
$\mu^{(l)}(\rho)_j^n$ and $U^{(l)}(\rho)^n_t$ as above, upon
substituting $c(\rho)$ and $u'_n/4$ with $c^{(l)}(\rho)$ and $u'_n/8$.
Since $U(\rho)^n_t\leq 4U^{(1)}(\rho)^n_t+4U^{(2)}(\rho)^n_t$,
it suffices to prove the result for each $U^{(l)}(\rho)^n_t$.

First, $\|\mu^{(1)}(\rho)^n_j\|^2\,1_{\{\|\mu^{(1)}(\rho)^n_i\|>u'_n/8\}}$
is smaller than $K\|\mu^{(1)}(\rho)^n_j\|^4/u'^2_n$, whereas (recalling
$\|\Wb(\rho)\|\leq K/\rho$) classical estimates yield
$\E\big(\|\mu^{(1)}(\rho)^n_j\|^4\big)\leq K\De_n(1+\De_n/\rho)$. Thus
the expectation of $U^{(1)}(\rho)^n_t$ is less than $K\De_n^{1/2-2\vpi'}
(1+\De_n/\rho)$, yielding the result for $U^{(1)}(\rho)^n_t$.

Secondly, we have $U^{(2)}(\rho)^n_t\leq\sum_{j=3}^{[t/k_n\De_n]}
\|\mu^{(2)}(\rho)^n_i\|^2$ and the first part of (\ref{P-11}) yields
$\E\big(\|\mu^{(2)}(\rho)^n_i\|^2\big)\leq K\phi_\rho\,\rdn$. Since
$\phi_\rho\to0$ as $\rho\to0$, the result for $U^{(1)}(\rho)^n_t$
follows.\qed

\subsection{A Basic Decomposition.}\label{ssec-2}

We start the proof of Theorem \ref{TR-1} by giving a decomposition of
$V(g)^n-V(g)$, with quite a few terms. It is
based on the key property $\wc'^n_i=c^n_i+\be^n_i$ and on the definition
\eqref{P-20} of $\al^n_i$ and $\be^n_i$. A simple
calculation shows that $\rdnn\,(V'(g)^n_t-V(g)_t)=\sum_{j=1}^5V^{n,j}_t$,
as soon as $t>k_n\De_n$, where
(the sums on components below always extend from $1$ to $d$):
$$\begin{array}{l}
V^{n,1}_t=\rdn\,\sdknt\big(g(\wc^n_i)-g(\wc'^n_i)\big)\\
V^{n,2}_t=\rdnn\,\sum\limits_{i=1}^{[t/\De_n]-k_n+1}
\int_{(i-1)\De_n}^{i\De_n}(g(c_i^n)-g(c_s))\,ds\\
V^{n,3}_t=\frac1{k_n\rdn}\,\sum\limits_{i=1}^{[t/\De_n]-k_n+1}
\sum\limits_{l,m} \partial_{lm}g(c^n_i)\,
\sum\limits_{u=0}^{k_n-1}\al^{n,lm}_{i+u}\\
V^{n,4}_t=\frac{\rdn}{k_n}\,\sum\limits_{i=1}^{[t/\De_n]-k_n+1}
\sum\limits_{l,m} \partial_{lm}g(c^n_i)\,
\sum\limits_{u=1}^{k_n-1}(c^{n,lm}_{i+u}-c^{n,lm}_i)
-\rdnn\,\int_{\De_n([t/\De_n]-k_n+1}^tg(c_s)\,ds\\
V^{n,5}_t=\rdn\,\sum\limits_{i=1}^{[t/\De_n]-k_n+1}
\big(g(c^n_i+\be^n_i)-g(c^n_i)
-\sum\limits_{l,m} \partial_{lm}g(c^n_i)\,\be^{n,lm}_i\big).
\end{array}$$

The leading term is $V^{n,3}$, the bias comes from the terms $V^{n,4}$
and $V^{n,5}$, and the first two terms are negligible, in the sense that
they satisfy
\bee\label{P-200}
j=1,2~~\Rightarrow~~V^{n,j}_t~\toop~0\qquad\text{for all $t>0$}.
\eee
We end this subsection with the proof of \eqref{P-200}.
\vsq

\nit The case $j=1$: \rm (\ref{R-8}) implies
$$|g(\wc^n_i)-g(\wc'^n_i)|\leq
K(1+\|\wc^n_i\|+\|\wc'^n_i\|)^{p-1}\,\|\wc^n_i-\wc'^n_i\|
\leq K(1+\|\wc'^n_i\|)^{p-1}\,\|\wc^n_i-\wc'^n_i\|
+K\|\wc^n_i-\wc'^n_i\|^p.$$
Recalling the last part of (\ref{P-21}), we deduce from (\ref{P-6}), from
the fact that $1-r\vpi-p(1-2\vpi)<\frac{(2-r)\vpi}{2q}$ for all
$q>1$ small enough, and from H\"older's
inequality, that $\E(|g(\wc^n_i)-g(\wc'^n_i)|)\leq
Ka_n\De_n^{(2p-r)\vpi+1-p}$. Therefore
$$\E\Big(\sup_{s\leq t}\,|V^{n,1}_s|\Big)\leq
Kta_n\De_n^{(2p-r)\vpi+1/2-p}$$
and (\ref{P-200}) for $j=1$ follows.
\vsq

\nit The case $j=2$: \rm Since $g$ is $C^2$ and $c_t$ is an It\^o
semimartingale with bounded characteristics, the convergence
$V^{n,2}\toucp0$ is well known: see for example the proof of
(5.3.24) in \cite{JP}, in which one replaces $\rho_{c_s}(f)$ by $g(c_s)$.

\subsection{The Leading Term $V^{n,3}$.}
Our aim here is to prove that
\bee\label{P-50}
V^{n,3}~\tolls~Z
\eee
(functional stable convergence in law), where $Z$ is the process
defined in Theorem \ref{TR-1}.

A change of order of summation allows us to rewrite $V^{n,3}$ as
$$V^{n,3}_t=\rdnn\sdt \sum_{l,m} w^{n,lm}_i\,\al^{n,lm}_i,\quad\text{where}~~
w^{n,lm}_i=\frac1{k_n}\sum_{j=(i-[t/\De_n]+k_n-1)^+}^{(i-1)\wedge(k_n-1)}
\partial_{lm} g(c_{i-j}^n).$$
Observe that $w^n_i$ and
$\al^n_i$ are measurable with respect to $\f^n_i$ and $\f^n_{i+1}$,
respectively, so by Theorem IX.7.28 of \cite{JS} (with $G=0$ and $Z=0$ in
the notation of that theorem) it suffices to prove the
following four convergences in
probability, for all $t>0$ and all component indices:
\bee\label{P-51}
\rdnn\sdknt w^{n,lm}_i\,\E(\al^{n,lm}_i\mid\f^n_i)~\toop~0
\eee
\bee\label{P-52}
\frac1{\De_n}\!\!\sdknt\!\! w^{n,jk}_i\,w^{n,lm}_i\,\E(\al^{n,jk}_i
\,\al^{n,lm}_i\mid\f^n_i)
\toop\int_0^t
\partial_{jk} g(c_s)\,\partial_{lm} g(c_s)\,
\big(c_s^{jl}c_s^{km}+c_s^{jm}c_s^{kl}\big)\,ds
\eee
\bee\label{P-53}
\frac1{\De_n^2}\sdknt\|w^n_i\|^4\,\E(\|\al^n_i\|^4\mid\f^n_i)
~\toop~0
\eee
\bee\label{P-54}
\rdnn\sdknt w^{n,lm}_i\,\E(\al^{n,lm}_i\,
\dd N\mid\f^n_i)~\toop~0,
\eee
where $N=W^j$ for some $j$, or is an arbitrary bounded
martingale, orthogonal to $W$.

For proving these properties, we pick a $\rho$ bigger than the upper
bound of the function $J$, so $\g^\rho$ becomes the trivial $\si$-field
and $\f^n_i=\f^{n,\rho}_i$ and $L(n,\rho)=\N$. In such a way, we can apply
all estimates
of the previous subsections with the conditioning $\si$-fields $\f^n_i$.
Therefore (\ref{P-21}) and the property $\|w^n_i\|\leq K$ readily
imply (\ref{P-51}) and (\ref{P-53}). In view of the form of $\al^n_i$,
a usual argument (see e.g. \cite{JP}) shows that in fact
$\E(\al^{n,lm}_i\,\dd N\mid\f^n_i)=0$ for all $N$ as above, hence
(\ref{P-54}) holds.

For (\ref{P-52}), by (\ref{P-22}) it suffices to prove that
$$\De_n\!\!\sdknt\!\! w^{n,jk}_i\,w^{n,lm}_i\,
(c^{n,jl}_ic^{n,km}_i+c^{n,jm}_ic^{n,kl}_i)
\toop\int_0^t\!\! \partial_{jk} g(c_s)\,\partial_{lm} g(c_s)\,
\big(c_s^{jl}c_s^{km}+c_s^{jm}c_s^{kl}\big)\,ds.$$
In view of the definition of $w^n_i$, for each $t$ we have
$w^{n,jk}_{i(n,t)}\to\partial_{jk} g(c_t)$ and $c^{n,jk}_{i(n,t)}\to
c^{jk}_t$ almost surely if $|i(n,t)\De_n-t|\leq k_n\De_n$ (recall that
$c$ is almost surely continuous at $t$, for any fixed $t$), and the above
convergence follows by the
dominated convergence theorem, thus ending the proof of (\ref{P-50}).

\subsection{The Term $V^{n,4}$.}\label{ssec-4}

In this subsection we prove that, for all $t$,
\bee\label{P-40}
V^{n,4}_t\toop \frac{\te}2\,\sum_{l,m}\int_0^t\partial_{lm}g(c_{s-})
\,dc_s^{lm}-\te\, g(c_t).
\eee

We call $V'^{n,4}_t$ and $V''^{n,4}_t$, respectively, the first
sum, and the last integral, in the definition of $V^{n,4}_t$.
Since $k_n\rdn\to\te$ and $c$ is a.s. continuous at $t$, it is obvious
that $V''^{n,4}_t$ converges almost surely to $-\te\, g(c_t)$, and
it remains to prove the convergence of $V'^{n,4}_t$ to the first
term in the right side of \eqref{P-40}.

We first observe that $c^{n}_{i+u}-c^{n}_i=\sum_{v=0}^{u-1}\De^n_{i+v}c$.
Then, upon changing the order of summation, we can rewrite $V'^{n,4}_t$ as
$$V^{n,4}_t=\sum_{i=1}^{[t/\De_n]-1}\sum_{l,m}w^{n,lm}_i\,\dd c^{lm},\quad
~~w^{n,lm}_i=\frac{\rdn}{k_n}\,\sum_{u=0\vee(i+k_n-1-[t/\De_n]}
^{i-1)\wedge(k_n-2)}(k_n-1-u)\partial_{lm}g(c^n_{i-u}).$$
In other words, recalling $k_n\rdn\leq K$ and $\|\partial g(c_s)\|\leq K$,
we see that
$$V^{n,4}_t=\sum_{l,m}\,\int_0^tH(n,t)^{lm}_s\,dc_s^{lm},$$
where $H(n,t)_s$ is a $d\times d$-dimensional predictable process, bounded
uniformly (in $n,s,\om$) and given on the set $[k_n\De_n,t-k_n\De_n]$ by
$$(i-1)\De_n<s\leq i\De_n~~\Rightarrow~~
H(n,t)_s^{lm}=\frac{\rdn}{k_n}\,\sum_{u=0}^{k_n-2}
(k_n-1-u)\partial_{lm}g(c^n_{i-u})$$
(its expression on $[0,k_n\De_n)$ and on $(t-k_n\De_n,t]$ is more
complicated, but not needed, apart from the fact that it is uniformly
bounded). Now, since $\sum_{u=0}^{k_n-2}(k_n-1-u)=k_n^2/2+$ O$(k_n)$ as
$n\to\infty$, we observe that $H(n,t)^{lm}_s$ converges to
$\frac{\te}2\,\partial_{lm}g(c_{s-})$ for all $s\in(0,t)$. Since
$c$ is a.s. continuous at $t$, we deduce from the dominated convergence
theorem for stochastic integrals that $V'^{n,4}_t$ indeed converges
in probability to the first term in the right side of \eqref{P-40}.

\subsection{The Term $V^{n,5}$.}\label{ssec-44}

The aim of this subsection is to prove the convergence
\bee\label{P-43}
V^{n,5}_t\toop A^2_t-2A^3_t+
\te\sum_{s\leq t}
\int_0^1\big(g(c_{s-}+w\De c_s)-g(c_{s-})-w
\sum_{l,m}\partial_{lm}g(c_{s-})\,\De c_s^{lm}\big)\,dw
\eee

We have $V^{n,5}_t=\sum_{i=1}^{[t/\De_n]-k_n+1}v^n_i$, where
$$v^n_i=\rdn\,\big(g(c^n_i+\be^n_i)-g(c^n_i)
-\sum\limits_{l,m} \partial_{lm}g(c^n_i)\,\be^{n,lm}_i\big).$$
We also set
\bee\label{P-430}
\begin{array}{c}
\Bal^n_i=\frac1{k_n\De_n}\sum_{u=0}^{k_n-1}\al^n_{i+u},\qquad
\Bbe^n_i=\be^n_i-\Bal^n_i=\frac1{k_n}\sum_{u=1}^{k_n-1}(c^n_{i+u}-c^n_i)\\
v'^n_i=\rdn\,\big(g(c^n_i+\Bbe^n_i)-g(c^n_i)
-\sum\limits_{l,m} \partial_{lm}g(c^n_i)\,\Bbe^{n,lm}_i\big),
\qquad v''^n_i=v^n_i-v'^n_i.
\end{array}
\eee

We take $\rho\in(0,1]$, and will eventually let it go to $0$.
With the sets $L(n,\rho)$ of (\ref{P-134}), we associate
$$\begin{array}{l}
L(n,\rho,t)=\{1,\cdots,[t/\De_n]-k_n+1\}\cap L(n,\rho)\\
\BL(n,\rho,t)=\{1,\cdots,[t/\De_n]-k_n+1\}\backslash L(n,\rho).
\end{array}$$
We split the sum giving $V^{n,5}_t$ into three terms:
\bee\label{P-431}
U^{n,\rho}_t=\sum_{i\in L(n,\rho,t)}v^n_i,\qquad
U'^{n,\rho}_t=\sum_{i\in \BL(n,\rho,t)}v'^n_i,\qquad
U''^{n,\rho}_t=\sum_{i\in \BL(n,\rho,t)}v''^n_i.
\eee
\vsq

\nib A) The processes $U^{n,\rho}$. \rm
A Taylor expansion and \eqref{R-8} give us
$$v^n_i=v(1)^n_i+v(2)^n_i+v(3)^n_i,~\text{where}~\left\{
\begin{array}{l}
v(1)^n_i=\frac{\rdn}2\,\sum_{j,k,l,m}\partial^2_{jk,lm} g(c^n_i)
\,\E(\be_i^{n,jk}\,\be_i^{n,lm}\mid\f^{n,\rho}_i)\\
v(2)^n_i=\frac{\rdn}2\,\sum_{j,k,l,m}\partial^2_{jk,lm} g(c^n_i)
\,\be_i^{n,jk}\,\be_i^{n,lm}-v(1)^n_i\\
|v(3)^n_i|\leq K{\rdn}\,(1+\|\be^n_i\|)^{p-3}\,\|\be^n_i\|^3.
\end{array}\right.$$
Therefore
\bee\label{P-432}
U^{n,\rho}=\sum_{j=1}^3U(j)^{n,\rho},\quad\text{where}\quad
U(j)^n_t=\sum_{i\in L(n,\rho,t)} v(j)^n_i.
\eee

On the one hand, and letting
$$w(\rho)^n_i=\sum_{j,k,l,m}\partial^2_{jk,lm} g(c^n_i)
\Big(\frac1{2k_n\rdn}\,(c^{n,jl}_ic^{n,km}_i+c^{n,jm}_ic^{n,kl}_i)
+\frac{k_n\rdn}6\,\Bc(\rho)_i^{n,jklm}\Big),$$
the c\`adl\`ag property of $c$ and $\Bc(\rho)$ and $k_n\rdn\to\te$ imply
$$W(\rho)^n_t:=\De_n\sdknt w(\rho)^n_i\toop U(1)^\rho_t:=
A^2_t+\frac{\te}6\sum_{j,k,l,m}\int_0^t\partial^2_{jk,lm}g(c_s)\,
\Bc(\rho)_s^{jklm}\,ds.$$
On the other hand, Lemma \ref{LP-8} yields
$[v(1)^n_i-\De_n w(\rho)^n_i|\leq K_\rho\De_n\,(\De_n^{1/4}+\eta(\rho)^n_i)$
when $i\in L(n,\rho)$, whereas $|w(\rho)^n_i|\leq K$ always. Therefore
$$\E\big(|U(1)^{n,\rho}_t-W(\rho)^n_t|\big)
\leq K_{\rho}\De_n\,\E\Big(\sdt(\rdn+\eta(\rho)^n_i)\Big)
+K\De_n\E\big(\#(\BL(n,\rho,t))\big).$$
Now, $\#(\BL(n,\rho,t))$ is not bigger than $(2k_n+1)N^\rho_t$, implying that
$\De_n\E(\#(L'(n,\rho,t)))\leq K_\rho\rdn$. Taking advantage of
\eqref{P-511}, we deduce that the above expectation goes to $0$ as
$n\to\infty$, and thus
\bee\label{P-42}
U(1)^{n,\rho}_t~\toop~U(1)^\rho_t.
\eee

Next, $v(2)^n_i$ is $\f^{n,\rho}_{i+k_n}$-measurable, with vanishing
$\f^{n,\rho}_i$-conditional expectation, and each set $\{i\in L(n,\rho)\}$
is $\f^{n,\rho}_0$-measurable. It follows that
$$\begin{array}{lll}
\E\big((U(2)^{n,\rho}_t)^2\big)&\leq &2k_n\,E\Big(
\sum_{i\in L(n,\rho,t)}\E\big(|v(2)^n_i|^2\mid\f^{n,\rho}_i\big)\Big)\\
&\leq& Kk_n\De_n\,E\Big(\sum_{i\in L(n,\rho,t)}\E\big(|\be^n_i|^4
\mid\f^{n,\rho}_i\big)\Big)~\leq~ Kt\phi_\rho+K_\rho t\rdn,
\end{array}$$
where we have applied (\ref{P-1401}) for the last inequality. Another
application of the same estimate gives us
$$\E\big(|U(3)^n_t|)~\leq~Kt\phi_\rho+K_\rho t\De_n^{1/4}.$$
These two results and the property $\phi_\rho\to0$ as $\rho\to0$
clearly imply
\bee\label{P-44}
\lim_{\rho\to0}\,\limsup_{n\to\infty}\,\E(|U(2)^{n,\rho}_t|
+|U(3)^{n,\rho}_t|)~=~0.
\eee
\vsq

\nib B) The processes $U'^{n,\rho}$. \rm We will use here the jump times
$S^\rho_1,S^\rho_2,\cdots$ of the Poisson process $N^\rho$,
and will restrict our attention to the set $\Om_{n,t,\rho}$ defined
before (\ref{P-101}), whose probability goes to $1$ as $n\to\infty$. On
this set, $\BL(n,\rho,t)$ is the collection of all integers $i$
which are between $[S^\rho_q/\De_n]-2k_n+2$ and $[S^\rho_q/\De_n]+1$, for
some $q$ between $1$ and $N^\rho_t$. Thus
\bee\label{P-441}
U'^{n,\rho}_t=\sum_{q=1}^{N^\rho_t}H(n,\rho,q),\quad\text{where}~~
H(n,\rho,q)=\sum_{i=[S^\rho_q/\De_n]-2k_n+1}
^{[S^\rho_q/\De_n]+1}v'^n_i.
\eee
The behavior of each $H(n,\rho,q)$ is a pathwise question. We fix $q$ and
set $S=S^\rho_q$ and $a_n=[S/\De_n]$, so $S>a_n\De_n$ because $S$ is not a
multiple of $\De_n$. For further reference we consider a case slightly
more general than strictly needed here. We have $c^n_i\to c_{S-}$ when
$a_n-6k_n+1\leq i\leq a_n+1$ and $c^n_i\to c_S$ when $a_n+2\leq i
\leq a_n+6k_n$, uniformly in $i$ (for each given outcome $\om$). Hence
\bee\label{P-49}
\Bbe^n_i-\frac{(k_n-a_n+i-2)^+\wedge(k_n-1)}{k_n}\,\De c_S\to0\quad
\text{uniformly in $i\in\{a_n-6k_n+2,\cdots,a_n+5k_n\}$.}
\eee
Thus, the following convergence holds, uniform in $i\in\{a_n-2k_n+1,\cdots,
a_n+1\}$:
$$\begin{array}{l}
\rdnn\,v'^n_i-\Big(g\big(c_{S-}+\frac{k_n-a_n+i-2}{k_n}
\,\De c_S\big)-g(c_{S-})\\ \hskip4cm
-\sum_{l,m}\partial_{lm}g(c_{S_-})\,\big(
c_{S-}^{jk}+\frac{k_n-a_n+i-2}{k_n}\,\De c_S^{lm}\big)\Big)\to0,
\end{array}$$
which implies
$$H(n,\rho,q)-\rdn\,\sum_{u=1}^{k_n-3}\Big(g\big(c_{S_q-}
+\frac u{k_n}\,\De c_{S_q}\big)-g(c_{S_q-})
-\sum_{l,m}\partial_{lm}g(c_{S_q-})\,\frac u{k_n}\,\De c^{lm}_{S_q}\Big)\to0$$
and by Riemann integration this yields
$$H(n,\rho,q)\to\te\int_0^1\big(g(c_{S_q-}+w\De c_{S_q})-g(c_{S_q-})-w
\sum_{l,m}\partial_{lm}g(c_{S_q-})\,\De c_S^{lm}\big)\,dw.$$
Henceforth, we have
\bee\label{P-45}
U'^{n,\rho}_t~\toop~U'^\rho_t:=\te\sum_{q=1}^{N^\rho_t}
\int_0^1\big(g(c_{S_q-}+w\De c_{S_q})-g(c_{S_q-})-w
\sum_{l,m}\partial_{lm}g(c_{S_q-})\,\De c_{S_q}^{lm}\big)\,dw.
\eee
\vsq

\nib C) The processes $U''^{n,\rho}$. \rm Since $|\Bbe^n_i|\leq K$ we deduce
from \eqref{R-8} that $|v''^n_i|\leq K\rdn\,(\|\Bal^n_i\|+\|\Bal^n_i\|^p)$.
\eqref{P-21} yields $\E(\|\Bal^n_i\|^q\mid\f^{n,\rho}_i)\leq K_q\De_n^{q/4}$
for all $q>0$. Therefore
$$\E\big(|U''^{n,\rho}_t|\big)~\leq~K\De_n^{3/4}\,\E(\#(\BL(,n,\rho,t)))
~\leq~K_\rho\De_n^{1/4},$$
by virtue of what precedes \eqref{P-42}. We then deduce
\bee\label{P-46}
U''^{n,\rho}_t~\toop~0.
\eee
\vsq

\nib D) Proof of \eqref{P-43}. \rm On the one hand, $V^{n,5}=U(1)^{n,\rho}
+U(2)^{n,\rho}+U(3)^{n,\rho}+U'^{n,\rho}+U''^{n,\rho}$; on the other hand,
the dominated convergence theorem (observe that $\Bc(\rho)_t\to\Wsi^2_t$
for all $t$) yields that $U(1)^\rho_t\toop A^2-\frac12\,A^3_t$ and
$$U'^\rho_t\toop \te\sum_{s\leq t}
\int_0^1\big(g(c_{s-}+w\De c_s)-g(c_{s-})-w
\sum_{l,m}\partial_{lm}g(c_{s-})\,\De c_s^{lm}\big)\,dw$$
as $\rho\to0$ (for the latter convergence, note that
$|g(x+y)-g(x)-\sum_{l,m}\partial_{lm}g(x)y^{lm}|
\leq K\|y\|^2$ when $x,y$ stay in a compact set).
Then the property (\ref{P-43}) follows from (\ref{P-42}), (\ref{P-44}),
(\ref{P-45}) and (\ref{P-46}).
\vsq

\nib E) Proof of Theorem \ref{TR-1}. \rm We are now ready to prove
Theorem \ref{TR-1}. Recall that $\rdnn\,(V(g)n_t-V(g))=\sum_{j=1}^5V^{n,j}$.
By virtue of \eqref{P-200}, \eqref{P-50}, \eqref{P-40}, \eqref{P-43},
it is enough to check that
$$\begin{array}{l}
A^1_t+A^3_t+A^4_t+A^5_t=\frac{\te}2\,\sum_{l,m}\int_0^t\partial_{lm}g(c_{s-})
\,dc_s^{lm}-\te\, g(c_t)\\
\qquad-2A^3_t+\te\sum_{s\leq t}
\int_0^1\big(g(c_{s-}+w\De c_s)-g(c_{s-})-w
\sum_{l,m}\partial_{lm}g(c_{s-})\,\De c_s^{lm}\big)\,dw.
\end{array}$$
To this aim, we observe that It\^o's formula gives us
$$g(c_t)=g(c_0)+\sum_{l,m}\int_0^t\partial_{lm}g(c_{s-})
\,dc_s^{lm}-\frac6{\te}A^3_t+
\sum_{s\leq t}\big(g(c_{s-}+\De c_s)-g(c_{s-})-
\sum_{l,m}\partial_{lm}g(c_{s-})\,\De c_s^{lm}\big),$$
so the desired equality is immediate (use also $\int_0^1w\,dw=\frac12$),
and the proof of Theorem \ref{TR-1} is complete.

\subsection{Proof of Theorem \ref{TR-2}.}\label{ssec-5}

The proof of Theorem \ref{TR-2} follows the same line as in Subsection
\ref{ssec-44}, and we begin with an auxiliary step.
\vsq

\nit Step 1) Replacing $\wc^n_i$ by $\wc'^n_i$. \rm The summands in
the definition (\ref{R-15}) of $A^{n,3}_t$ are $R(\wc^n_i,
\wc^n_{i+k_n})$, where $R(x,y)=\sum_{j,k,l,m}\partial_{jk,lm}^2g(x)
(y^{jk}-x^{jk})(y^{lm}-x^{lm})$, and we set
$$A'^{n,3}_t=-\frac{\rdn}8\,
\sum_{i=1}^{[t/\De_n]-2k_n+1}R(\wc'^n_i,\wc'^n_{i+k_n}).$$
We prove here that
\bee\label{P-70}
A^{n,3}_t-A'^{n,3}_t~\toop~0
\eee
for all $t$, and this is done as in to the step $j=1$ in Subsection
\ref{ssec-2}. The function $R$ is $C^1$ on $\R_+^2$ with
$\|\partial^jR(x,y)\|\leq K(1+\|x\|+\|y\|)^{p-j}$ for $j=0,1$, by
(\ref{R-8}). Thus
$$\begin{array}{lll}
|R(\wc^n_i,\wc^n_{i+k_n})-R(\wc'^n_i,\wc'^n_{i+k_n})|&\leq&
K(1+\|\wc'^n_i\|+\|\wc'^n_{i+k_n})\|)^{p-1}(\|\wc^n_i-\wc'^n_i\|
+\|\wc^n_{i+k_n}-\wc'^n_{i+k_n}\|)\\&&\qquad
+K\|\wc^n_i-\wc'^n_i\|^p+K\|\wc^n_{i+k_n}-\wc'^n_{i+k_n}\|^p.
\end{array}$$
Then, exactly as in the case afore-mentioned, we conclude (\ref{P-70}),
and it remains to prove that, for all $t$, we have
$$A'^{n,3}_t~\toop~-\frac12\,A^{2}_t+A^3_t+A'^{4}_t.$$
\vsq

\nit Step 2) \rm From now on we use the same notation as in Subsection
\ref{ssec-44}, although they denote different variables or processes.
For any $\rho\in(0,1]$ we have $A'^{n,3}=U^{n,\rho}+U'^{n,\rho}+U''^{n,\rho}$,
as defined in \eqref{P-431}, but with
$$\begin{array}{l}
v^n_i=-\frac{\rdn}8\,R(c^n_i+\be^n_i,c^n_{i+k_n}+\be^n_{i+k_n})\\
v'^n_i=-\frac{\rdn}8\,R(c^n_i+\Bbe^n_i,c^n_{i+k_n}+\Bbe^n_{i+k_n}),\quad
v''^n_i=v^n_i-v'^n_i.
\end{array}$$

Recalling $\ga^n_i$ in \eqref{P-20}, the decomposition (\ref{P-432}) holds with
$$\begin{array}{l}
v(1)^n_i=-\frac{\rdn}8\,\sum_{j,l,k,m}\partial^2_{jl,km}\,g(c^n_i)\,
\,\E\big(\ga_i^{n,jk}\,\ga_i^{n,lm}\mid\f^{n,\rho}_i\big)\\
v(2)^n_i=-\frac{\rdn}8\,\sum_{j,l,k,m}\partial^2_{jl,km}\,g(c^n_i)
\,\ga_i^{n,jk}\,\ga_i^{n,lm}-v(1)^n_i\\
v(3)^n_i=v^n_i-v(1)^n_i-v(2)^n_i.
\end{array}$$
Use $\wc'^n_i-c^n_i=\be^n_i$ and (\ref{R-8}) and
a Taylor expansion to check that
$$|v(3)^n_i|\leq K\rdn\,\|\ga^n_i\|^2\,\|\be^n_i\|\,
(1+\|\be^n_i\|)^{p-3}.$$
We also have $|v(2)^n_i|\leq K\rdn\,\|\ga^n_i\|^2$, hence (\ref{P-1401})
and (\ref{P-26}) yield
$$\E(|v(3)^n_i|\mid\g^\rho)+
\E(|v(2)^n_i|^2\mid\g^\rho)\leq K\De_n\,\Big(\phi_\rho+\De_n^{1/4}
+\frac{\De_n}{\rho^p}\Big),$$
and thus (\ref{P-44}) holds here as well, by the same argument. Moreover,
(\ref{P-26}) again
yields (\ref{P-42}), with now
$$U^\rho_t=-\sum_{j,k,l,m}\int_0^t\partial^2_{jk,lm}\,g(c_s)\,
\Big(\frac{\te}{12}\,\Bc(\rho)^{jklm}_s+\frac1{4\te}\,
(c_s^{jl}c_s^{km}+c_s^{jm}c_s^{kl})\Big)\,ds.$$
This goes to $A^3_t-\frac12\,A^2_t$ as $\rho\to0$.

Another application of \eqref{R-8} gives us
$$|v''^n_i|\leq K\rdn\,\big(1+\|\ga^n_i\|^2\big)
\big(\|\Bal^n_i\|+\|\Bal^n_{i+k_n}\|+
\|\Bal^n_i\|^p+\|\Bal^n_{i+k_n}\|^p\big).$$
Then another application of \eqref{P-21}, \eqref{P-1401} and
\eqref{P-26} yields $\E(|v''^n_i|\mid\g^\rho)\leq K\De_n^{3/4}$ and we
conclude \eqref{P-46} as previously. We are thus left to prove that
\bee\label{P-71}
\rho>0~\Rightarrow~U'^{n,\rho}_t\toop U'^\rho_t,\quad
\text{with, as $\rho\to0$,}~~U'^\rho_t\toop A_t'^4.
\eee
\vst

\nit Step 3) \rm On the set $\Om_{n,t,\rho}$ we have \eqref{P-441}
and we study $H(n,\rho,q)$, in the same way as before, on the set
$\Om_{n,t,\rho}$. We fix $q$ and set $S=S_q$ and $a_n=[S/\De_n]$.
We then apply \eqref{P-49} and also $c^n_i\to c_{S-}$ or $c^n_i\to c_S$,
according to whether $a_n-2k_n+1\leq i\leq a_n+1$ or $a_n+2\leq i\leq a_n+k_n$,
to obtain $v'^n_i-\Bv'^n_i\to0$, uniformly in $i$
between $a_n-2k_n+1$ and $a_n+1$, where
$$\Bv'^n_i=\left\{\begin{array}{l}
0\hskip5cm \text{if}~a_n-2k_n+1\leq i\leq a_n-2k_n+2\\
-\frac{(2k_n-a_n+i-2)^2\rdn}{8k^2_n}\sum_{j,k,l,m}\partial^2_{jk,lm}\,g(c_{S-})
\,\De c_S^{jk}\,\De c_S^{lm}\\
\hskip5cm~~\text{if}~a_n-2k_n+3\leq i\leq a_n-k_n+1\\
\frac{(a_n-i+2)^2\rdn}{8k^2_n}\sum_{j,k,l,m}\partial^2_{jk,lm}\,g\Big(
c_{S-}+\frac{k_n-a_n+i+2}{k_n}\,\De c_S\Big)
\,\De c_S^{jk}\,\De c_S^{lm}\\
\hskip5cm~~\text{if}~a_n-k_n+2\leq i\leq a_n+1.
\end{array}\right.$$
We then deduce, by Riemann integration, that
$$H(n,\rho,q)\to-\frac{\te}8\!\sum_{j,k,l,m}\int_0^1\!\!
\big(\partial^2_{jk,lm}\,g(c_{S_q-})+
\partial^2_{jk,lm}\,g(c_{S_q-}+(1-w)\De c_{S_q})\big)\,w^2
\,\De c_{S_q}^{jk}\,\De c_{S_q}^{lm}\,dw,$$
which is $\te G'(c_{S_q-},\De c_{S_q})$, hence
the first part of (\ref{P-71}), with $U'^\rho_t=
\te\sum_{q=1}^{N^\rho_t}G'(c_{S^\rho_q-},\De c_{S^\rho_q})$. The
second part of (\ref{P-71}) follows from the dominated
convergence theorem, and the proof of Theorem \ref{TR-2} is complete.

\subsection{Proof of Theorem \ref{TR-3}.}

The proof is once more somewhat similar to the proof of Subsection
\ref{ssec-44}, although the way we replace $\wc^n_i$ by $\wc'^n_i$ and
further by $\Bal^n_i+\be^n_i$ is different.
\vsq

\nib A) Preliminaries. \rm The $j$th summand in (\ref{R-22}) involves
several estimators $\wc^n_i$, spanning the time interval
$((j-3)k_n\De_n,(j+2)k_n\De_n]$. It is thus convenient
to replace the sets $L(n,\rho)$,
$L(n,\rho,t)$ and $\BL(n,\rho,t)$, for $\rho,t>0$, by the following ones:
$$\begin{array}{l}
L'(n,\rho)=\{j=3,4,\cdots:~N^\rho_{(j+2)k_n\De_n}-N^\rho_{(j-3)k_n\De_n}
=0\}\\
L'(n,\rho,t)=\{3,\cdots,[t/k_n\De_n]-3\}\cap L'(n,\rho)\\
\BL'(n,\rho,t)=\{3,\cdots,[t/k_n\De_n]-3\}\cap(\N\backslash L'(n,\rho)).
\end{array}$$

For any $\rho\in(0,1]$ we write
$\va(F)^n_t=\va^{n,\rho}_t+\Bva^{n,\rho}_t$, where
$$\begin{array}{c}
v^n_j=F(\wc^n_{(j-3)k_n+1},\de^n_j\wc)\,1_{\{\|\de^n_{j-1}\wc\|\vee
\|\de^n_{j+1}\wc\|\vee u'_n<\|\de^n_j\wc\|\}}\\[1.6mm]
\va^{n,\rho}_t=\sum_{j\in L'(n,\rho,t)}v^n_j,\qquad
\Bva^{n,\rho}_t=\sum_{j\in \BL'(n,\rho,t)}v^n_j.
\end{array}$$
We also set
$$\begin{array}{ll}
\de^n_j\wc'=\wc'^n_{jk_n+1}-\wc'^n_{(j-2)k_n+1},\qquad&
\de^n_j\be=\be^n_{jk_n+1}-\be^n_{(j-2)k_n+1}\\
w^n_j=\sum_{m=-3}^2\|\wc^n_{(j+m)k_n+1}-\wc'^n_{(j+m)k_n+1}\|,\quad&
w'^n_j=(1+\|\wc'^n_{(j-3)k_n+1}\|)^{p-1}\,(1+\|\de^n_j\wc\|)^2.
\end{array}$$

(\ref{P-6}) and the last part of (\ref{P-21}) yield
\bee\label{P-80}
q\geq1~~\Rightarrow~~\E\big((w^n_j)^q) \leq
K_q\,\De_n^{(2q-r)\vpi+1-q},\qquad \E((w'^n_i)^q)\leq K_q.
\eee
Observe that $\de^n_j\wc'$ is analogous to $\ga^n_i$, with a doubled
time lag, so it satisfies a version of (\ref{P-26}) and, for $q\geq2$,
we have
\bee\label{P-81}
i\in L'(n,\rho)~\Rightarrow~
\E\big(\|\de^n_j\wc'\|^q\mid\f^{n,\rho}_{(j-2)k_n+1})\big|\leq K_q
\big(\rdn\,\phi_\rho+\De_n^{q/4}+\frac{\De_n^{q/2}}{\rho^q}\big).
\eee
\vsq

\nib B) The processes $\va^{n,\rho}$. \rm (\ref{R-14}) yields
$$|v^n_j|\leq K(1+\|\wc^n_{(j-3)k_n+1}\|)^{p-2}\,
\|\de^n_j\wc\|^2\,1_{\{\|\de^n_j\wc\|>u'_n\}}+K\|\de^n_j\wc\|^p.$$
Thus a (tedious) computation shows that, with the notation
$$a^n_j=(1+\|\wc'^n_{(j-3)k_n+1}\|)^{p-2}\,
\|\de^n_j\wc'\|^2\,1_{\{\|\de^n_j\wc'\|>u'_n/2\}},\quad a'^n_j=
w'^n_j\Big(w^n_i+(w^n_i)^p+\frac{(w^n_i)^v}{u_n'^v}\Big),$$
with $v>0$ arbitrary,
we have $|v^n_j|\leq K(a^n_j+\|\de^n_j\wc'\|^p+a'^n_j)$ (with $K$
depending on $v$). Therefore
we have $|\va^{n,\rho}_t|\leq K(B^{n,\rho}_t+C^{n,\rho}_t+D^n_t)$,
where
$$B^{n,\rho}_t=\sum_{j\in L'(n,\rho,t)}a^n_j,\qquad
C^{n,\rho}_t=\sum_{j\in L'(n,\rho,t)}\|\de^n_j\wc'\|^p,\qquad
D^n_t=\sum_{j=3}^{[t/k_n\De_n]}a'^n_j.$$

First, (\ref{P-80}) and H\'older's inequality give us $\E(a'^n_j)
\leq K_{q,v}\De_n^{l(q,v)}$ for any $q>1$ and $v>0$, where (recalling
(\ref{R-9}) and (\ref{R-20}) for $\vpi$ and $\vpi'$) we have set
$l(q,v)=\frac{1-r\vpi}q-\big(p(1-2\vpi)\vee v(1-2\vpi+\vpi')\big)$.
Upon choosing $v$ small enough and $q$ close enough to $1$, and in view of
(\ref{R-9}), we see that $l(q,v)>\frac12$, thus implying
\bee\label{P-84}
\E(D^n_t)~\to~0.
\eee
Next, we deduce from (\ref{P-81}) that
$$\E\big(C^{n,\rho}_t\big)\leq K
\E\Big(\E\Big(
\sum_{i\in L'(n,\rho,t)}\|\de^n_j\wc'\|^p\mid\g^\rho\Big)\Big)
\leq Kt\Big(\phi_\rho+\De_n^{p/4}+\frac{\De_n^{p/2}}{\rho^p}\Big),$$
and thus, since $p\geq3$,
\bee\label{P-85}
\lim_{\rho\to0}\,\limsup_{n\to\infty}\,\E(|C^{n,\rho}_t|)~=~0.
\eee

The analysis of $B_t^{n,\rho}$ is more complicated. We have
$\de^n_j\wc'=z^n_j+z'^n_j$, where
$$z^n_j=\Bal^n_{jk_n+1}-\Bal^n_{(j-2)k_n+1},\quad
z'^n_j=\frac1{k_n}\,\sum_{m=1}^{k_n}(c^n_{jk_n+m}-c^n_{(j-2)k_n+m})$$
(recall \eqref{P-430} for $\Bal^n_i$), hence
$$a^n_j\leq 4(1+\|\wc'^n_{(j-3)k_n+1}\|)^{p-2}\,\Big(
\|z^n_j\|^2\,1_{\{\|z^n_i\|>u'_n/4\}}
+\|z'^n_j\|^2\,1_{\{\|z'^n_i\|>u'_n/4\}}\Big).$$
It easily follows that for all $A>1$,
\bee\label{P-86}
B^{n,\rho}_t~\leq~16\,B^{n,\rho,1}_t+4A^{p-2}\,B^{n,\rho,2}_t
+\frac{2^p}A\,B^{n,\rho,3}_t,
\eee
where
$$\begin{array}{c}
B^{n,\rho,m}_t=\sum_{j\in L'(n,\rho,t)}a(m)^n_j,\qquad
a(1)^n_j=(1+\|\wc'^n_{(j-3)k_n+1}\|)^{p-2}\,\frac{\|z^n_j\|^3}{u'_n}\\
a(2)^n_j=\|z'^n_j\|^2\,1_{\{\|z'^n_i\|>u'_n/4\}},\quad
a(3)^n_j=\|\wc'^n_{(j-3)k_n+1}\|^{p-1}\,\|z'^n_j\|^2.
\end{array}$$

On the one hand, (\ref{P-21}) and H\"older's inequality yield
$\E(a(1)^n_j\mid\g^{\rho})\leq K\De_n^{3/4-\vpi'}$ and, since
$\vpi'<\frac14$, we deduce
\bee\label{P-87}
\E\big(B^{n,\rho,1}_t\big)~\to~0.
\eee
On the other hand, observe that $z'^n_j=\mu(\rho)^n_j$, with
the notation (\ref{P-30}), and as soon as $j\in L'(n,\rho)$, so
Lemma \ref{LP-2} gives us
\bee\label{P-88}
\lim_{\rho\to0}~\limsup_{n\to\infty}~\E\big(B^{n,\rho,2}_t\big)~=~0.
\eee
Finally, (\ref{P-11}) shows that $\E(\|z'^n_j\|^q|\mid
\f^{n,\rho}_{(j-2)k_n+1})\leq
K_{q,\rho}\rdn$ for all $q\geq2$ and $j\in L'(n,\rho)$, whereas
$\wc'^n_{(j-3)k_n+1}$ is $\f^n_{(j-2)k_n+1}$-measurable, so (\ref{P-11}),
(\ref{P-21}) and successive conditioning yield
$\E(a(3)^n_j|\mid\g^\rho)\leq K_{q,\rho}\rdn$.
Then, again as for (\ref{P-87}), one obtains
\bee\label{P-89}
\E\big(B^{n,\rho,3}_t\big)~\leq~K_\rho\,t.
\eee

At this stage, we gather (\ref{P-84})--(\ref{P-89}) and obtain, by
letting first $n\to\infty$, then $\rho\to0$, then $A\to\infty$, that
\bee\label{P-90}
\lim_{\rho\to0}~\limsup_{n\to\infty}~\E\big(|\va^{n,\rho}_t|\big)~=~0.
\eee
\vsq

\nib C) The processes $\Bva^{n,\rho}$. \rm With the previous notation
$S^\rho_j$ and $N^\rho_t$, and on the set $\Om_{n,\rho,t}$,
we have
\bee\label{P-91}
\Bva^{n,\rho}_t=\sum_{m=1}^{N^\rho_t}~
\sum_{j=-2}^2v^n_{[S^\rho_m/k_n\De_n]+j}\,.
\eee
This is a finite sum (bounded in $n$ for each $\om$). Letting $S=
S^\rho_m$ for $m$ and $\rho$ fixed and $w_n=\frac S{k_n\De_n}-
\big[\frac S{k_n\De_n}\big]$, we know that for any given $j\in\Z$ the
variable $\wc^n_{([S/k_n\De_n]+j)k_n+1}$ converge in probability to $c_{S-}$
if $j<0$ and to $c_S$ if $j>0$, whereas for $j=0$ we have
$\wc^n_{[S/k_n\De_n]k_n+1}-w_nc_S-(1-w_n)c_S\toop0$. This in turn implies
$$\begin{array}{c}
j<0~\text{or}~j>2~\Rightarrow~\de^n_{[S/k_n\De_n]+j}\wc\toop0\\
\de^n_{[S/k_n\De_n]}\wc-(1-w_n)\De c_S\toop0,\quad
\de^n_{[S/k_n\De_n]+1}\wc\toop\De c_S,\quad
\de^n_{[S/k_n\De_n]+2}\wc-w_n\De c_S\toop0.
\end{array}$$
By virtue of the
definition of $v^n_j$, and since $u'_n\to0$ and also since
$w_n$ is almost surely in $(0,1)$ and $F$ is continuous and $F(x,0)=0$,
one readily deduces that
$$v^n_{[S/k_n\De_n]+j}~\toop~\left\{\begin{array}{ll}
F(c_{S-},\De c_S)~~&\text{if}~j=1\\ 0&\text{if}~j\neq1.
\end{array}\right.$$

Coming back to (\ref{P-91}), we deduce that
\bee\label{P-93}
\Bva^{n,\rho}_t~\toop~\Bva^\rho_t:=\sum_{m=1}^{N^\rho_t}
F(c_{S^\rho_m-},\De c_{S^\rho_m}).
\eee
In view of (\ref{R-14}), an application of the dominated convergence
theorem gives $\Bva^\rho_t\to \va(F)_t$. Then (\ref{R-24}) follows
from $\va(F)^n_t=\va^{n,\rho}_t+\Bva^{n,\rho}_t$ and (\ref{P-90})
and (\ref{P-93}), and the proof of Theorem \ref{TR-3} is complete.


\begin{thebibliography}{9}

\bibitem{APPS} Alvarez, A., Panloup, P., Pontier, M. and Savy, N. (2010).
Estimation of the instantaneous volatility. {\em Statistical Inference for
Stochastic Processes} {\bf 15}, 27-59.

\bibitem{CDG} Cl\'ement, E., Delattre, S. and Gloter, A. (2012).
An infinite dimensional convolution theorem with
applications to the efficient estimation of the integrated volatility.
Preprint.


\bibitem{JS}
Jacod, J. and Shiryaev, A.N. (2003).
\textit{Limit Theorems for Stochastic Processes}, 2nd ed.
Springer-Verlag, Berlin.

\bibitem{JP} Jacod, J. and Protter, P. (2012).
\textit{Discretization of Processes}, Springer-Verlag, Berlin.

\bibitem{JR} Jacod, J. and Rosenbaum, M. (2012). Quarticity and other
functionals of volatility: efficient estimation.

\bibitem{V10} Vetter, M. (2010). Limit theorems for bipower variation
of semimartingales. {\em Stochastic Processes and their Applications}
{\bf 120}, 22-38.

\bibitem{V1} Vetter, M. (2011). Estimation of integrated volatility of
volatility with applications to a goodness-of-fit testing. Preprint.

\end{thebibliography}
\end{document}